\documentclass[11pt]{article}

  \newcommand{\lab}[1]{\label{#1}}                % hides labels

\usepackage[latin1]{inputenc}
\usepackage{amssymb,amsmath,amsthm}
\usepackage{mathrsfs}
\newcommand{\no}{\noindent}

\theoremstyle{plain}
\newtheorem{thm}{Theorem}[section]
\newtheorem{lem}[thm]{Lemma}
\newtheorem{cor}[thm]{Corollary}
\newtheorem{prop}[thm]{Proposition}
\theoremstyle{definition}
\theoremstyle{remark}

\newenvironment{pf}{\begin{proof}[\bf Proof]}{\end{proof}}

\newcommand{\real}{\ensuremath {\mathbb R} }

\newcommand{\st}{ : } % Try this: the other makes too much space.

\newcommand{\ex} {{\bf E}}
\newcommand{\pr} {{\bf P }}

\newcommand{\mbf}[1] {\text{\boldmath$#1$}}
\newcommand{\remove}[1] {}

\newcommand{\cD} {{\mathcal D}}

\newcommand{\cG} {{\mathcal G}}
\newcommand{\cH} {{\mathcal H}}
\newcommand{\cP} {{\mathcal P}}
\newcommand{\cQ} {{\mathcal Q}}
\newcommand{\hD} {{\widehat{\mathcal D}}}

\newcommand{\Gnm}{{\cG(n,m)}}
\newcommand{\dico}{{\cG_{1,1}(n,m)}}
\newcommand{\dicok}{{\cG_k(n,m)}}
\newcommand{\Pnm}{{\cP_{1,1}(n,m)}}
\newcommand{\Pknm}{{\cP_k(n,m)}}
\newcommand{\Psnm}{{\cP'_{1,1}(n,m)}}
\newcommand{\Hnm}{{\cH(n,m)}}
\newcommand{\Qnm}{{\cQ(n,m)}}

\newcommand{\dv} {{\vec{d}}}
\newcommand{\dpv} {{\vec{d^+}}}
\newcommand{\dmv} {{\vec{d^-}}}
\newcommand{\ba}{{\mbf a}}

\newcommand{\Gd}{{\cG(\dv)}}
\newcommand{\Pd}{{\cP(\dv)}}
\newcommand{\Psd}{{\cP'(\dv)}}
\newcommand{\Qd}{{\cQ(\dv)}}
\newcommand{\Qa}{{\cQ(\ba)}}
\newcommand{\Hd}{{\cH(\dv)}}
\newcommand{\Hsd}{{\cH'(\dv)}}
\newcommand{\Ha}{{\cH(\ba)}}
\newcommand{\Hsa}{{\cH'(\ba)}}

\newcommand{\good}{{\Gamma}}
\newcommand{\Bad}{\ensuremath{\mathcal B}}

\newcommand{\tpo}{{\mathrm{TPo}}}
\newcommand{\tpol}{{\mathrm{TPo}(\lambda)}}
\newcommand{\bindis}{{\rm Bin}}
\newcommand{\distrib}{\overset{\mathrm d}{\sim}}

\newcommand{\cN} {{\mathcal N}}
\newcommand{\Nscr}{\cN}

\newcommand{\hatdico}{{\hat\cG_{1,1}(n,m,n_1)}}

\newcommand\la{\lambda}
\newcommand{\eps}{{\epsilon}}

\def\tb{{\bf t}}

\newcommand{\be}{\begin{equation}}
\newcommand{\ee}{\end{equation}}
\newcommand{\bea}{\begin{eqnarray}}
\newcommand{\eea}{\end{eqnarray}}
\newcommand{\bean}{\begin{eqnarray*}}
\newcommand{\eean}{\end{eqnarray*}}
 \newcommand\eqn[1]{(\ref{#1})}
\newcommand{\bel}[1]{\be\lab{#1}}

\title{Asymptotic enumeration of strongly connected digraphs by vertices and edges}
\author{Xavier P\' erez-Gim\' enez\thanks{Partially supported by the Province of Ontario under the Post-Doctoral Fellowship (PDF) Program.} \ and
Nicholas Wormald\thanks{Supported by the Canada Research Chairs Program,   NSERC Discovery Grants Program, and  CRM, Bellatera, Spain.} \\
{\small Department of Combinatorics and Optimization} \\
{\small University of Waterloo}\\
{\small Waterloo ON, Canada}\\
\smallskip
{\small {\tt xperez@uwaterloo.ca}, {\tt nwormald@uwaterloo.ca }}
 }
\date{}

\begin{document}
\maketitle
\begin{abstract}
We derive an asymptotic formula for the number of strongly connected digraphs with $n$ vertices and $m$ arcs (directed edges), valid for $m-n\to\infty$ as $n\to \infty$ provided $m=O(n\log n)$. This fills the gap between Wright's results which apply to $m=n+O(1)$, and the long-known threshold for $m$, above which a random digraph with $n$ vertices and $m$ arcs is likely to be strongly connected.
\end{abstract}
\section{Introduction} \lab{s:intro}
One of the most fundamental properties of a directed graph (digraph), and possibly the most useful for communication networks, is that of being {\em strongly connected}, that is, possessing directed paths both ways between every pair of vertices. It was long ago shown  by Moon and Moser~\cite{MM} that almost all of the $2^{n^2}$ digraphs with $n$ vertices are strongly connected due to having paths of length 2 between each pair of vertices.  So, from the asymptotic enumeration perspective, a more interesting problem is the enumeration of strongly connected digraphs on $n$ vertices with   $m$  arcs (i.e.\ directed edges). In this paper, all digraphs are labelled. Our results cover only \emph{simple} digraphs (i.e.\ digraphs with no multiple arcs), but unless otherwise stated we allow digraphs to have loops. We also give results for digraphs in which loops are forbidden, which we refer to as \emph{loop-free} digraphs.

Pal{\' a}sti~\cite{P} determined the threshold  of strong connectivity, as follows. Let $\alpha$ be fixed and define $m(\alpha,n)= \lfloor n\log n+ \alpha  n\rfloor$. Then, for a random directed graph having $n$ vertices and $m$ arcs, so that each of the ${n^2 \choose N}$ possible choices is equiprobable, the probability that the digraph is strongly connected tends to $\exp(-2e^{-\alpha})$ as $n\to\infty$. Multiplying this probability by ${n^2 \choose N}$ consequently gives an asymptotic formula for the number   $S(n,m)$ of strongly connected digraphs with $n$ vertices and $m$ arcs, for such $m$. This also easily implies that  $S(n,m)\sim {n^2 \choose N}$ if $m=m(\alpha_n,n)$ with $\alpha_n\to \infty$.
 On the other hand, Wright~\cite{Wright} obtained recurrences for  the exact value of $S(n,m)$   when $m=n+O(1)$. (We must require $m\ge n$ to avoid the   failure to be strongly connected for trivial reasons.) In this paper, we fill   the entire gap  between these results, deriving an asymptotic formula  for $S(n,m)$, valid for $m-n\to\infty$ as $n\to \infty$ provided $m=O(n\log n)$. Our main result is as follows.
\begin{thm}\lab{t:main} Uniformly for $m=O(n\log n)$ and $m-n\to \infty$, the number of strongly connected
digraphs with $n$ vertices and $m$ arcs is asymptotic to
\begin{equation}
\frac{(m-1)! (e^\lambda-1)^{2n}}{2\pi(1+\lambda-c) \lambda^{2m}}\exp
(-\lambda^2/2 )
\frac{e^\lambda(e^\lambda-1-\lambda)^2}
{(e^{2\lambda}-e^\lambda-\lambda)(e^\lambda-1)},
\lab{eq:main}
\end{equation}
where $c=m/n>1$ and $\lambda$ is determined by the equation $c=\lambda e^\lambda/(e^\lambda-1)$.
\end{thm}
\paragraph{Note.}
In particular, if $c\to1$ then the   expression~\eqref{eq:main} simplifies asymptotically to
\begin{equation}
\frac{(m-1)! (e^\lambda-1)^{2n}}{6\pi \lambda^{2m}},
\lab{eq:main2}
\end{equation}
whilst if $c\to\infty$ then~\eqref{eq:main} is asymptotic to
\begin{equation}
\frac{(m-1)! (e^\lambda-1)^{2n}}{2\pi \lambda^{2m}}\exp
(-\lambda^2/2 ).
\lab{eq:main3}
\end{equation}
 
Our result has counterparts for undirected graphs. An asymptotic formula for the number of connected graphs with $n$  vertices and $m$ edges was   given  for all $m$ such that $m-n\to\infty$ as $n\to\infty$ by Bender, Canfield and McKay~\cite{BCM}. This improved the range of $m$ for which earlier estimates were found, and also the bounds on the error term.  %
A simpler approach to the same problem was given in~\cite{PW2}. This begins by counting connected graphs with no end-vertices, and then consider the number of ways to attach a forest. One of the ways used there to count connected 2-cores was to count connected kernels, which have no vertices of degree 2, and insert vertices of degree 2 into their edges, and another way was based on eliminating isolated cycles by inversion. In the present paper, for the case $m=O(n)$ we use this first of these two alternatives. This has some advantage in providing direct information on properties of the kernel, such as was used in~\cite{KW} for studying long cycles in the supercritical random graph. In a similar way, we can study the analogous structure for a digraph, which we call its heart. For $m/n\to \infty$ we use a rather different approach to show that random digraphs with all in- and out-degrees at least 1 are strongly connected with high probability.

Our argument requires a formula for the number of digraphs with all in- and outdegrees at least $1$ and given number of arcs, which we obtain using the method for counting graphs with given minimum degree developed by Pittel and the second author in~\cite{PW1}. 
\begin{thm}\lab{t:dicoreenum}
Uniformly for $m=O(n\log n)$,
%
%\begin{equation}\lab{didegseq}
$$
|\dico| \sim \frac{m! (e^\lambda-1)^{2n}}{2\pi n c(1+\lambda-c) \lambda^{2m}}\exp
(-\lambda^2/2 )
$$
%\end{equation}
%
where $c=m/n$ and $\lambda$ is determined by $c=\lambda e^\lambda/(e^\lambda-1).$
\end{thm}
\noindent
Using the same method, we also extend this to digraphs with all outdegrees at least $k^+$ and all indegrees at least $k^-$.
\remove{Along the way, we obtain a formula for the number of digraphs with all outdegrees at least $k^+$, indegrees at least $k^-$, and given number of arcs,  using the method for counting graphs with given minimum degree developed by Pittel and the second author in~\cite{PW1}.}
\begin{thm} \lab{t:kdicoreenum}
Fix positive integers $k^+$ and $k^-$. Uniformly for $m=O(n\log n)$, $m-k^+n\to+\infty$ and $m-k^-n\to+\infty$, the number of digraphs on $n$ vertices, $m$ arcs, outdegrees at least $k^+$ and indegrees at least $k^-$ is asymptotic to
\[
\frac{(m-1)! (f_{k^-}(\lambda^-)f_{k^+}(\lambda^+))^n}{2\pi \sqrt{(1+\eta^+-c)(1+\eta^--c)} \lambda^{2m}}\exp
(-\lambda^-\lambda^+/2),
\]
where 
$c=m/n,
$ 
$$
f_k(\lambda)=\sum_{i\ge\max\{k,0\}}\lambda^i/i!,
$$  
$\lambda^+$ and $\lambda^-$ are the unique positive roots of
$$
c = \lambda^+ f_{k^+-1}(\lambda^+)/f_{k^+}(\lambda^+),\quad
c = \lambda^- f_{k^--1}(\lambda^-)/f_{k^-}(\lambda^-)
$$
respectively, and
$$
\eta^+ = (\lambda^+)^2 f_{k^+-2}(\lambda^+)/f_{k^+-1}(\lambda^+),\quad \eta^- = (\lambda^-)^2 f_{k^--2}(\lambda^-)/f_{k^--1}(\lambda^-).
$$
\end{thm}

The results stated so far refer to digraphs that are allowed to have loops but not multiple arcs. In Section~\ref{s:loopfree} we extend these results to the case when loops are forbidden, and obtain the following analogues of Theorems~\ref{t:main} and~\ref{t:kdicoreenum}.
\begin{thm}\lab{t:mainloopfree} Uniformly for $m=O(n\log n)$ and $m-n\to+\infty$, the number of strongly connected
loop-free digraphs with $n$ vertices and $m$ arcs is asymptotic to
\[
\frac{(m-1)! (e^\lambda-1)^{2n}}{2\pi(1+\lambda-c) \lambda^{2m}}\exp
(-c(1- e^{-\lambda})^2-\lambda^2/2 )
\frac{e^\lambda(e^\lambda-1-\lambda)^2}
{(e^{2\lambda}-e^\lambda-\lambda)(e^\lambda-1)},
\]
where $c$ and $\lambda$ are as in is Theorem~\ref{t:main}.
\end{thm}
\noindent Note that, for Theorem~\ref{t:mainloopfree}, the only effect of forbidding loops was   to introduce the extra factor $\exp (-c(1- e^{-\lambda})^2)$.
\begin{thm} \lab{t:kdicoreenumloopfree}
Fix positive integers $k^+$ and $k^-$, and recall  the notation of Theorem~\ref{t:kdicoreenum}. Uniformly for $m=O(n\log n)$, $m-k^+n\to+\infty$ and $m-k^-n\to+\infty$, the number of loop-free digraphs on $n$ vertices, $m$ arcs, outdegree at least $k^+$ and indegree at least $k^-$ is asymptotic to
\[
\frac{(m-1)! (f_{k^-}(\lambda^-)f_{k^+}(\lambda^+))^n}{2\pi \sqrt{(1+\eta^+-c)(1+\eta^--c)} \lambda^{2m}}\exp
(-c-\lambda^-\lambda^+/2).
\]
\end{thm}
\noindent For Theorem~\ref{t:kdicoreenumloopfree}, forbidding loops just gave only the factor $e^{-c}$.

Cooper and Frieze~\cite[Theorem 3(vi)]{CF} obtained a significant result relevant to this problem, in the form of the asymptotic probability that a random digraph with given degree sequence is strongly connected, under certain assumptions on the degree sequence. It would be rather straightforward to combine this with our Theorem~\ref{t:dicoreenum}, along with properties of   degree sequences which we use in our paper,  to deduce an asymptotic formula for  $S(n,m)$ when   $m/n>1 $ is bounded away from 1 and is bounded.  For completeness, we derive this case of the   formula in a different way, following the same approach as we use for the case $m/n\to 1$, which we consider in Section~\ref{s:vsscd}.

Boris Pittel~\cite{Boris} has independently investigated the second approach of~\cite{PW2} mentioned above. Applying it to this problem in the loop-free case, he has simultaneously obtained a formula similar to that in Theorem~\ref{t:mainloopfree} under the restriction that $m=O(n)$, but also including an explicit error estimate.
\section{Basics and notation}\lab{s:basics}
\subsection{Truncated Poisson distribution}
We consider a discrete probability distribution that will be used many times in the argument. Given $\lambda>0$ and a nonnegative integer $k$, we say that a random variable (r.v.)\ $Y$ has a {\em $k$-truncated Poisson distribution of parameter $\lambda$} (or simply $Y\distrib\tpo_k(\lambda)$) if
\[
\pr(Y=i)=\begin{cases}
\displaystyle \frac{\lambda^i}{f_k(\lambda)\,i!} & \text{if } i\ge k,
\\
0 & \text{if } 0\le i<k,
\end{cases}
\]
where $f_k(\lambda)=\sum_{i\ge k}\lambda^i/i!$.  For later convenience we also define $f_k(\lambda)=e^\lambda$ for integer $k< 0$.

We first give a rough tail bound for a random variable   $Y\distrib\tpo_k(\lambda)$ for constant $k$ but $\lambda$ possibly depending on $n$. Consider constants $A>B>e$, and let $p$ be a constant nonnegative integer. Then for $j\ge\max\{A e\lambda,k\}$ we have
\begin{align*}
\ex([Y]_p\,1_{Y\ge j}) &=
\frac{1}{f_k(\lambda)}\left(\sum_{j\le i<j+p}[i]_p\frac{\lambda^i}{i!} +\sum_{i\ge j+p}[i]_p\frac{\lambda^i}{i!}\right)
\le \frac{1}{f_k(\lambda)}\left(p(j+p)^p\frac{\lambda^j}{j!} +\lambda^p\sum_{i\ge j}\frac{\lambda^i}{i!}\right)
\notag\\
&= O\left( \frac{j^p+\lambda^p}{f_k(\lambda)} (e\lambda/j)^j \right) = O\left(B^{-j}\right).
\end{align*}
In particular,
\begin{equation}
\pr(Y\ge j) = O\left(B^{-j}\right),
\quad
\ex(Y\,1_{Y\ge j}) = O\left(B^{-j}\right)
\quad\text{and}\quad
\ex([Y]_2\,1_{Y\ge j}) = O\left(B^{-j}\right).
\lab{eq:tailbounds}
\end{equation}
(We use $[x]_k$ to denote the falling factorial $x(x-1)\cdots (x-k+1)$  throughout this paper.)

Our main use   of the $\tpo_k(\lambda)$ distribution is to allow us to make computations on the multinomial distribution truncated from below. The following lemma establishes a connection between these distributions, and will be used throughout the paper often without an explicit mention. (See for example~\cite[Section~2]{CW} for a proof of this lemma.)
\begin{lem}\lab{l:multinomial}
Distribute $M\ge kN$ distinguishable balls randomly into $N$ distinguishable bins u.a.r.\ subject to the condition that each bin receives at least $k\ge1$ balls. Let $Y_i$ be the numbers of balls in bin $i$. Then the joint distribution of $Y_1,\ldots,Y_N$ is the same as that of $N$ independent copies of $\tpo_k(\lambda)$ for arbitrary $\lambda>0$ conditional upon $Y_1+\cdots+Y_N=M$.
\end{lem}
It is easy to see that a variable $Y\distrib\tpo_k(\lambda)$   has   $\ex Y=c$ given by
\begin{equation}
c = \frac{\lambda f_{k-1}(\lambda)}{f_k(\lambda)}.
\lab{eq:ladef}
\end{equation}
Henceforth, given $c>k$, we assume that $\lambda$ is set equal to the unique (by~\cite[Lemma~1]{PW1}) positive root of this equation.
We also define
\begin{equation}
\eta = \frac{\lambda^2 f_{k-2}(\lambda)}{f_{k-1}(\lambda)}.
\lab{eq:etadef}
\end{equation}
Elementary computations show that, for such choice of $\lambda$ and $\eta$, we have  $\ex(Y(Y-1)) = \eta c$. More properties of the $\tpo_k(\lambda)$ distribution are given in~\cite{PW1}. It is easy to check that $0<\lambda\le c$ in all cases.
From~\cite[Theorem~4(a)]{PW1} we have the following.
\begin{lem}\lab{l:sum}
Let $M=O(N\log N)$ be integer such that $r:=M-kN\to\infty$ and put $c=M/N$. Let $Y_1,\ldots,Y_N$ be i.i.d.\ random variables with $\tpo_k(\lambda)$ distribution, for fixed $k$, where $\la$ is determined from $c$ in~\eqn{eq:ladef}, and define $\eta$ as in~\eqn{eq:etadef}.   Then, as $N\to\infty$,
\[
\pr(Y_1+\cdots+Y_N=M)\sim \frac{1}{\sqrt{2\pi Nc(1+\eta-c)}} = \Theta\left(1/\sqrt{r}\right).
\]
\end{lem}
Throughout the paper, we mostly focus our attention to the case $k=1$ and simply refer to the $\tpo_1(\lambda)$ distribution as $\tpol$ or simply {\em truncated Poisson}. In this particular case, \eqref{eq:ladef} can be rewritten as
\bel{eq:ladef1}
c=\frac{\lambda e^\lambda}{e^\lambda-1},
\ee
and moreover we have $\eta=\lambda$.

On several occasions we use   Chernoff bounds for a binomially distributed ${\rm Bin}(n,p)$ random variable $X$ in the common form
\bel{chern}
\pr(|X-np|>a)<2e^{-2a^2/n},
\ee
or the variation more useful when $p$ is small:
\bel{chernsmallp}
\pr(|X-np|>a)<2e^{-a^2/3np} \quad \mbox{for} \quad   a\le np
\ee
(from Molloy~\cite{Mo}; see also Alon and Spencer~\cite[Theorems~A.1.11 and~A.1.13]{AS08}).

We close this subsection with some rather technical lemmas on independent variables with $\tpo_k(\lambda)$ distribution.
\begin{lem}\lab{l:concentration}
Let $Y_1,\ldots,Y_N$ be independent r.v.s\ with $\tpo_k(\lambda)$ distribution, for fixed $k$ and for $0<\lambda\le\log N$. Put $C=\ex Y_1$. Then for any $t\ge\sqrt N\log^2N$ we have
\[
\pr\left(\Big|\sum_{i=1}^N Y_i - CN\Big|>t\right) = O\left( e^{-(t^2/8N)^{1/3}} \right),
\]
asymptotically as $N\to\infty$.
\end{lem}
\begin{pf}
Let $Y_{\max}=\max_i \{Y_i\}$. Setting $\Delta=(t^2/8N)^{1/3}$, we have $\pr(Y_{\max}>\Delta) \le N \pr(Y_1>\Delta) = O(e^{-\Delta})$ by~\eqref{eq:tailbounds} and since $\Delta=\Omega(\log^{4/3}N)$. Now define
\[
W_i = Y_i-C \quad\text{and}\quad W^*_i = W_i \,1_{Y_i\le \Delta},
\]
and again from~\eqref{eq:tailbounds} deduce
\begin{equation}
|\ex W^*_i| = |-\ex(W_i \,1_{Y_i>\Delta})|
\le \ex(Y_i \,1_{Y_i>\Delta}) = O\left(e^{-\Delta}\right).
\lab{eq:EWi}
\end{equation}
Moreover, we have $-C \le W^*_i \le \Delta-C$, and then $|W^*_i-\ex W^*_i|<\Delta$, so by the Azuma-Hoeffding inequality
\begin{equation}
\pr\bigg( \Big|\sum_{i=1}^n (W^*_i - \ex W^*_i) \Big| \ge t/2 \bigg)
\le 2 \exp\left(\frac{-t^2}{8\Delta^2N}\right) = 2 e^{-\Delta}.
\lab{eq:Azuma2}
\end{equation}
\begin{align*}
\pr\left(\Big|\sum_{i=1}^N Y_i-CN\Big|>t\right) &\le
\pr(Y_{\max}>\Delta) + 
\pr\bigg( \Big|\sum_{i=1}^N W^*_i \Big| > t \bigg)
\\
&\le O\left( e^{-\Delta} \right) +\pr\bigg( \Big|\sum_{i=1}^N (W^*_i - \ex W^*_i) \Big| > t - \Big|\sum_{i=1}^N \ex W^*_i \Big| \bigg)
\\
&\le O\left( e^{-\Delta} \right) + \pr\bigg( \Big|\sum_{i=1}^N (W^*_i - \ex W^*_i) \Big| > t/2 \bigg)
\\
&= O\left( e^{-\Delta} \right),
\end{align*}
where we used~\eqref{eq:Azuma2} and the fact that $|\sum_{i=1}^N \ex W^*_i|<t/2$, which follows from~\eqref{eq:EWi}.
\end{pf}
The following result was essentially shown in~\cite{PW1}.
\begin{lem}\lab{l:variance}
Let $Y_1,\ldots,Y_N$ be independent r.v.s\ with $\tpo_k(\lambda)$ distribution, for fixed $k$ and for $0<\lambda\le\log N$. Put $C=\ex(Y_1(Y_1-1))$. Then
\[
\pr\left(\Big|\sum_{i=1}^N Y_i(Y_i-1) - CN\Big|>4N^{1/2}\log^8N\right) = O\big(\exp(-\log^3 N)),
\]
asymptotically as $N\to\infty$.
\end{lem}
\begin{pf}
The statement in the lemma comes directly from equation~(33) in~\cite{PW1}, considering  (16), (22), (28), (29) and Lemmas~1 and 2 of that paper. See also the proof of Lemma~\ref{l:concentration} which uses the same method in full detail.
\end{pf}
We will use the following for $k=1,2$.
\begin{lem}\lab{l:binsfromballs}
Let $k\ge1$ be an integer, and let $Y_1,\ldots,Y_N$ be independent $\tpo_k(\lambda)$ r.v.s. Consider $N$ bins, place $Y_i$ balls in bin $i$ ($i=1,\ldots,N$), and then select each ball independently with probability $ q\le1/2$ where $Nq \ge \log^2 N$. Then the number $X$ of bins containing at least one selected ball satisfies
\[
\pr(|X-\ex X|>\sqrt{\ex X}\log N) = e^{-\Omega(\log^2N)}
\]
asymptotically as $N\to\infty$, and moreover 
$$\ex X/n>kq(1-(k-1)/4 + (2^{-k}/k)\pr(Y_1\ge k+1)).$$
\end{lem}
\begin{pf}
Let $q'$ be the probability that a bin contains at least one selected ball. We have
\begin{align*}
1-q' &< (1-q)^k\pr(Y=k)+(1-q)^{k+1}\pr(Y_i\ge k+1)
\\
&= (1-q)^k-q(1-q)^k\pr(Y_i\ge k+1).
\end{align*}
Using the elementary bound $(1-q)^k\le1-kq+\binom{k}{2}q^2$ and the fact that $q\le1/2$, we obtain

\begin{align}
q' &> kq(1-(k-1)q/2) - q(1-q)^k\pr(Y_i\ge k+1)
\notag\\
&\ge kq(1-(k-1)/4 + (2^{-k}/k)\pr(Y_i\ge k+1)),
\lab{eq:qprime}
\end{align}
and trivially $q'\ge q$ in any case.
Since $X\distrib\bindis(N,q')$, it follows by~\eqn{chernsmallp} that
\begin{equation}
\pr\big(|X-Nq'|>\sqrt{Nq'}\log N\big) < 2e^{-\log^2N/3}.\qedhere
\label{eq:chernoff}
\end{equation}
\end{pf}
%
%

%%%%%%%%%%%%%%%%%%%%%%%%%%%%%%%%%%%%%
\subsection{Probability spaces of digraphs and degree sequences}
Let $\Gnm$ be the set of digraphs on $n$ labelled vertices and $m$ arcs. In our definition of digraph we allow loops but not multiple arcs. It is a simple matter to adjust our arguments for  loop-free digraphs (see Section~\ref{s:loopfree}).
% query Add the formula!
  For a given digraph in $\Gnm$, let $\dpv=(d^+_1,\ldots,d^+_n)$ and $\dmv=(d^-_1,\ldots,d^-_n)$ denote respectively the sequences of out- and indegrees of the vertices.
The degree of vertex $i$ is defined to be the tuple $d_i=(d^+_i,d^-_i)$, so the joint in- and outdegree sequences can be represented by $\dv=(d_1,\ldots,d_n)$.
For feasibility, it is necessary that
\begin{equation}\lab{eSum}
\sum_{i=1}^nd^+_i=\sum_{i=1}^nd^-_i=m.
\end{equation}
Let $c=m/n$ and assume that $c>1$ throughout the article, though $m$ and hence $c$ are functions of $n$. Let $\dico$ be the set of digraphs in $\Gnm$ such that $d^+_i,d^-_i\ge1$ for all $i\in\{1,\ldots,n\}$. (Note that this is a necessary condition for strong connectedness when $n>1$.) Elements of $\dico$ we call $(1,1)$-\emph{dicores} or simply \emph{dicores}. We also write $\Gnm$ and $\dico$ to denote the corresponding uniform probability spaces. We define $r=m-n=(c-1)n$ and assume $r\to\infty$. 
We distinguish three subcases: very sparse, with $r=o(n)$ or equivalently $c\to1$; 
 moderately sparse, with $r=\Theta(n)$; and a denser case, with $c\to\infty$ but 
  $c=O(\log n)$. (All logarithms are natural unless otherwise specified.)

Let $\cD$ be the set of sequences $\dv=(d_1,\ldots,d_n)$, with $d_i=(d^+_i,d^-_i)$ for $i\in\{1,\ldots,n\}$, where the $2n$ entries $d^+_i$ and $d^-_i$ are positive integers. Let $\hD$ be the subset of sequences in $\cD$ satisfying the total sum conditions~\eqref{eSum}. Note that $\hD$ coincides with the set of all possible degree sequences of dicores in $\dico$. Given any $\dv\in\hD$, let $\Gd$ denote the set (and also the corresponding uniform probability space) of digraphs with degree sequence $\dv$. Also consider the usual directed pairing model $\Pd$, defined as follows. Take $n$ bins, where the $i$-th bin contains points of two types, namely $d^+_i$ \emph{out-points} and $d^-_i$ \emph{in-points}, and consider a random matching of the $m$ out-points with the $m$ in-points. Each element in $\Pd$ corresponds to a multidigraph in the obvious way, and the restriction to \emph{simple} digraphs (i.e.\ with no multiple arcs) generated this way is uniform.

In order to study the distribution of degree sequences of $\dico$, it will prove useful to turn the sets $\cD$ and $\hD$ into suitable probability spaces, as follows. Random degree sequences $\dv\in\cD$ are chosen by taking the $2n$ entries $d^+_i$ and $d^-_i$ as independent copies of $\tpol$. Let $\Sigma$ be the event in $\cD$ that~\eqref{eSum} holds, and define $\hD$ to be the corresponding conditional probability space. Moreover, let $\Pnm$ be the probability space of random pairings in $\Pd$ where the degree sequence $\dv$ is drawn from the distribution of $\hD$ defined above. Each pairing in $\Pnm$ corresponds to a multidigraph, and as will become apparent later the restriction of $\Pnm$ to simple digraphs generates elements of $\dico$ uniformly.

We also need the   notation  $d^+_{\max} =\max \{d^+_i  \st  1\le i\le n\}$ and $d^-_{\max} =\max \{d^-_i \st 1\le i\le n\}$.

%%%%%%%%%%%%%%%%%%%%%%%%%%%%%%%%%%%%%%%%%%%%%%%%%%%%%%%%%%%%%%%%%%%%%%%%%%%%%%%%%%%%
\section{Asymptotic enumeration of dicores}
Here we prove Theorems~\ref{t:dicoreenum} and~\ref{t:kdicoreenum} by adapting the main argument of~\cite{PW1}. Before that, we need some lemmata. The following result is an immediate consequence of Theorem~4.6 in~\cite{M84} by McKay (we just need to use the standard interpretation of digraphs with loops as bipartite graphs).
\begin{lem}[McKay]\lab{l:PSimple}
Let $\dv\in\hD$ be a sequence of degrees and suppose that $d^+_{\max},d^-_{\max}\le\Delta$ for some $\Delta=o(m^{1/4})$. Then the probability that a random element of $\Pd$ has no multiple arcs is
\[
\exp\left(-\frac{1}{2m^2}\sum_{i,j=1}^nd^+_i(d^+_i-1)d^-_j(d^-_j-1) + O\left(\frac{\Delta^4}{m}\right)\right),
\]
uniformly for all $\dv$.
\end{lem}
The following technical result estimates the probability that a degree sequence in $\cD$ satisfies~\eqref{eSum}, and averages the probability that a random pairing is simple over any subset of degree sequences with that property. Here $\la$ and $c$ are defined as in Theorem~\ref{t:dicoreenum}.
\begin{lem}\lab{l:useful}\hspace{0cm}
Assume that $m-n\to\infty$ and $m=O(n\log n)$.
\begin{itemize}
\item[(a)] $\displaystyle \pr_\cD(\Sigma) \sim \frac{1}{2\pi nc(1+\lambda-c)}=\Theta(1/(m-n))$.
\end{itemize}
Moreover, if $S$ is the event that a random pairing in $\Pd$ or $\Pnm$ is simple, then
\begin{itemize}
\item[(b)] $\displaystyle \pr_{\Pnm}(S) = \ex_\hD \big( \pr_\Pd(S)\big) \sim e^{-\lambda^2/2}$;
\item[(c)] for any r.v.\ $X$ on $\hD$ satisfying $|X|\le x$ for some fixed constant $x\in\real$,
\[
\ex_\hD \big( \pr_\Pd(S) \cdot X \big) =
(1 + o(1)) \, e^{-\lambda^2/2}\, \ex_\hD X + O\left(e^{-\log^3n}\right).
\]
\end{itemize}
\end{lem}
\begin{pf}
{From} Lemma~\ref{l:sum}, the independent events $\sum_id^+_i=m$ and $\sum_id^-_i=m$ each have probability $(1+o(1))/\sqrt{2\pi nc(1+\lambda-c)}$, which gives (a).
Note that (b) follows from (c) by setting $X=1$, since the bound on $m$ implies $\la=O(\log n)$, so it only remains to prove (c). For this, we follow  the proof of~\cite[Theorem~4(b)]{PW1}  almost exactly.

We require some definitions. Let $F=F(\dv)=\pr_\Pd(S)$ and
\[
\tilde F = \exp\left(-\frac12 D^+D^-\right),
\]
where
\[
D^+ = \frac{1}{m}\sum_{i=1}^nd^+_i(d^+_i-1)\qquad\text{and}\qquad
D^- = \frac{1}{m}\sum_{j=1}^n d^-_j(d^-_j-1).
\]
We set $\Delta=\log^3n$, and let $\Bad_1$ denote the `bad' event that $d^+_{\max}>\Delta$ or $d^-_{\max}>\Delta$. {From}~\eqref{eq:tailbounds} we obtain $\pr_\cD(\Bad_1) \le 2n \pr(Y>\Delta)=O(nA^{-\Delta})$. Then, we use the result from (a) to deduce that $\pr_\hD(\Bad_1) \le \pr_\cD(\Bad_1)/\pr_\cD(\Sigma) = O(n^2c(1+\lambda-c)A^{-\Delta}) = O\big(e^{-\log^3n}\big)$.

In view of Lemma~\ref{l:PSimple} and bearing in mind that $0\le F,\tilde F\le1$ and $|X|\le x$, we can write
\begin{align}
\ex_\hD(FX) &=
\ex_\hD(FX\,1_{\Bad_1}) + \ex_\hD(FX\,1_{\overline{\Bad_1}})
\notag\\
&= O(\pr_\hD(\Bad_1)) + (1+O(\Delta^4/m)) \ex_\hD(\tilde FX\,1_{\overline{\Bad_1}})
\notag\\
&= O\big(e^{-\log^3n}\big) + (1+O(\Delta^4/m)) \ex_\hD(\tilde FX\,1_{\overline{\Bad_1}}).
\lab{eq:EFBad1}
\end{align}
Simple computations show that $\ex D^+=\ex D^-=\lambda$ (with $D^+$ and $D^-$ independent). Set $t=8n^{-1/2}\log^9n$, and define $\Bad_2$ to be the `bad' event that $|D^+D^-/2-\lambda^2/2|>t$. Whenever $\Bad_2$ does not hold, we have $\tilde F = \exp(-\lambda^2/2 + O(t)) = (1 + O(t))\exp(-\lambda^2/2)$, so
\begin{align}
\ex_\hD(\tilde FX\,1_{\overline{\Bad_1}}) &=
\ex_\hD(\tilde FX\,1_{\overline{\Bad_1}\wedge\Bad_2}) + \ex_\hD(\tilde FX\,1_{\overline{\Bad_1}\wedge\overline{\Bad_2}})
\notag\\
&=
O(\pr_\hD(\Bad_2))
+ (1+O(t))\, e^{-\lambda^2/2}\, \ex_\hD(X)
\lab{eq:EFBad2}.
\end{align}
It only remains to bound $\pr_\hD(\Bad_2)$. Set $s=t/(2\log n)=4n^{-1/2}\log^8n$, and note that if $|D^+-\lambda|\le s$ and $|D^--\lambda|\le s$ then
\[
|D^+D^-/2-\lambda^2/2| \le \frac{1}{2}(|D^+-\lambda||D^--\lambda|+\lambda|D^+-\lambda|+\lambda|D^--\lambda|) \le \frac{s^2+2s\log n}{2} \le t.
\]
Therefore, by Lemma~\ref{l:variance},
\begin{equation}
\pr_\hD(\Bad_2) \le \pr_\hD(|D^+-\lambda|>s) + \pr_\hD(|D^--\lambda|>s) = O\big(e^{-\log^3n}\big).
\lab{eq:S1S2}
\end{equation}
Part (c) in the statement follows by combining~\eqref{eq:EFBad1}, \eqref{eq:EFBad2} and~\eqref{eq:S1S2}.
\end{pf}
Now we are in good shape to prove the theorem.%
\begin{proof}[\bf Proof of Theorem~\ref{t:dicoreenum}]
Observe that $|\Pd|=m!$, and that each simple digraph with degree sequence $\dv$ comes from exactly $\prod_{i=1}^n d^+_i!d^-_i!$ different pairings in $\Pd$. Thus
\[
|\Gd| = \frac{m! \pr_\Pd(S)} {\prod_{i=1}^n d^+_i!d^-_i!},
\]
where $S$ denotes the event that a random pairing in $\Pd$ has no multiple arcs.
Define
\[
Q = \sum_{\dv\in\hD} \prod_{i=1}^n\frac{1}{d^+_i!d^-_i!}
= \frac{(e^\lambda-1)^{2n}}{\lambda^{2m}} \pr_{\cD}\left(\Sigma\right).
\]
Therefore, summing over all degree sequences, we can write
\begin{align*}
|\dico| &=  \sum_{\dv\in\hD}
\frac{m! \pr_\Pd(S)} {\prod_{i=1}^n d^+_i!d^-_i!}
\\
&= m! Q
\ex_{\hD}\left(\pr_\Pd(S)\right)
\\
&= m! \frac{(e^\lambda-1)^{2n}}{\lambda^{2m}}
\ex_{\Pnm}(S) \pr_{\cD}\left(\Sigma\right)
\\
&\sim \frac{m! (e^\lambda-1)^{2n}}{2\pi nc(1+\lambda-c) \lambda^{2m}}
\exp(-\lambda^2/2),
\end{align*}
where we used Lemma~\ref{l:useful}.
\end{proof}
In addition, the computations in the proof of Theorem~\ref{t:dicoreenum} give the following.
\begin{cor}\lab{c:PtoG}
The elements in $\dico$ can be uniformly generated by restricting the probability space $\Pnm$ to simple pairings and considering the corresponding digraph. 
\end{cor}
\begin{pf}
A dicore $G$ in $\dico$ with degree sequence $\dv$ comes from exactly $\prod_{i=1}^n d^+_i!d^-_i!$ different pairings. Each of these pairings must be simple and has probability
\begin{equation}\lab{eq:PtoG}
\frac{\lambda^{2m}/(e^\lambda-1)^{2n}}{m!\prod_{i=1}^n d^+_i!d^-_i!}
 \big(\pr_\Pnm(S)\big)^{-1}
\end{equation}
in the space $\Pnm$ conditional upon the event $S$ of being simple. The product of~\eqref{eq:PtoG} times $\prod_{i=1}^n d^+_i!d^-_i!$ does not depend on the particular $\dv$, and therefore the distribution of $G$ when generated from simple pairings is uniform.  
\end{pf}
Finally, we can extend the concept of dicore defined in Section~\ref{s:basics} as follows. Given $k=(k^+,k^-)$ where $k^+$ and $k_-$ are positive integer constants, a $k$-\emph{dicore} is an element of $\Gnm$ with a degree sequence satisfying $d^+_i\ge k^+$ and $d^-_i\ge k^-$, for all $i\in\{1,\ldots,n\}$. Let $\dicok$ denote both the set of $k$-dicores and the corresponding uniform probability space.

In order to study the degree sequences of $\dicok$, we need some definitions. Let $\lambda^+$ and $\eta^+$ (resp., $\lambda^-$ and $\eta^-$) be obtained from~\eqref{eq:ladef} and~\eqref{eq:etadef} after replacing $k$, $\lambda$ and $\eta$ by $k^+$, $\lambda^+$ and $\eta^+$ (resp., by $k^-$, $\lambda^-$ and $\eta^-$). Define the set of degree sequences $\cD_k$ analogously to $\cD$, with the extra condition that $d^+_i\ge k^+$ and $d^-_i\ge k^-$, for all $i\in\{1,\ldots,n\}$, and similarly let $\hD_k$ be the subset of sequences in $\cD_k$ satisfying~\eqref{eSum}. Moreover, we endow $\cD_k$ with a probability distribution by selecting the $d^+_i$ and the $d^-_i$ independently according to the $\tpo_{k^+}(\lambda^+)$ and the $\tpo_{k^-}(\lambda^-)$ distributions, respectively. The $\hD_k$ space is simply $\cD_k$ conditional upon~\eqref{eSum}. Furthermore, we define $\Pknm$ as we did for $\Pnm$ but randomising the degree sequence $\dv$ according to the distribution of $\hD_k$ defined above.

Now we are in good shape to extend the argument in the proof of Theorem~\ref{t:dicoreenum} to general $k$-dicores.
\begin{proof}[\bf Proof of Theorem~\ref{t:kdicoreenum}]
The proof is straightforward by going along the  same steps as the proof of Theorem~\ref{t:dicoreenum}, but replacing $\cD$, $\hD$ and $\Pnm$ by  $\cD_k$, $\hD_k$ and $\Pknm$, and considering the distributions $\tpo_{k^+}(\lambda^+)$ or $\tpo_{k^-}(\lambda^-)$ instead of $\tpol$ when appropriate. The key part is extending Lemma~\ref{l:useful} to the new setting, which is also straightforward. The extended statement is as follows. Assume that $m-k^+n\to\infty$, $m-k^-n\to\infty$ and $m=O(n\log n)$. Then
\begin{itemize}
\item[(a)] $\displaystyle \pr_{\cD_k}(\Sigma) \sim \frac{1}{2\pi nc\sqrt{(1+\eta^+-c)(1+\eta^--c)}}$.
\end{itemize}
Moreover, if $S$ is the event that a random pairing in $\Pd$ or $\Pknm$ is simple, then
\begin{itemize}
\item[(b)] $\displaystyle \pr_{\Pknm}(S) = \ex_{\hD_k} \big( \pr_\Pd(S)\big) \sim e^{-\lambda^+\lambda^-/2}$;
\item[(c)] for any r.v.\ $X$ on $\hD_k$ satisfying $|X|\le x$ for some fixed constant $x\in\real$,
\[
\ex_{\hD_k} \big( \pr_\Pd(S) \cdot X \big) =
(1 + o(1)) \, e^{-\lambda^+\lambda^-/2}\, \ex_{\hD_k} X + O\left(e^{-\log^3n}\right).
\qedhere
\]
\end{itemize}
\end{proof}
%
%-----------------------------------------
%-----------------------------------------
%
\section{Moderately sparse case: $c$ bounded}\lab{s:sscd}

In  this section we will prove Theorem~\ref{t:main} for the case that $c=m/n$ is bounded and also bounded away from $1$.

A {\em sink-set} in a digraph $G$ is a non-empty proper subset $S$ of vertices such that the out-set of $S$ is a subset of $S$. That is, no arc goes from $S$ to $V(G)\setminus S$. A set of vertices is a {\em source-set} if its complement is a sink-set.
A sink-set in a digraph with minimum outdegree at least $1$ is {\em plain} if its vertices all have outdegree exactly $1$, and is otherwise {\em complex}. Plain and complex source-sets are defined analogously by replacing outdegree by indegree.
Observe that a digraph $G$ is strongly connected iff it has no sink-set (and equivalently no source-set). We use the term \emph{s-set} to denote sets of vertices which are a sink-set or a source-set.

We first show that a.a.s.\ any complex s-set of $\dico$ must contain more than $m/2$ arcs. Therefore, the strong connectedness of $\dico$ can be characterised in terms of plain s-sets.
\begin{prop}\lab{p:complex} Suppose that $c=m/n$ is bounded and bounded away from $1$. A digraph in
$\dico$ a.a.s.\ has no complex s-set containing at most $m/2$ arcs.
\end{prop}
\begin{pf}
It is presumably possible to analyse $\dico$ or $\Pnm$ directly to achieve the desired result, by an expectation argument similar to that commonly used for connectivity of graphs. However, the expectation itself seems to be difficult to analyse. Instead we introduce another probability space, by partitioning according to the indegree sequence and to the multiset of outdegrees. More precisely, we will consider slices of $\Pnm$ with indegree sequence $\dmv$ and outdegree sequence being a permutation of $\dpv$, for each $\dv\in\cD$.

One could argue by partitioning according to the joint values of
$\dmv$ and $\dpv$, but certain nasty combinations of in- and outdegrees, in which the vertices of outdegree 1 all have large indegree, are likely to cause trouble, and rather ad-hoc arguments may be required to bound the troublesome cases (see e.g.~the approach in~\cite{CF}). It is conceivable that allowing permutations of the outdegree sequence instead helps to explain a little more of the structure of the typical digraph in $\Gnm$.

To facilitate calculations of probabilities, for each $\dv\in\hD$ we introduce a probability space, $\Psd$, which is similar to common models (called pairing or configuration models) for random graphs or digraphs with given degree sequence. Consider two sets of points $A=\{a_1, \ldots, a_m\}$ and $B=\{b_1, \ldots , b_m\}$, with $A$ partitioned into nonempty sets (which we call {\em bins}) $A_i$, $i=1,\ldots ,n$ (corresponding to the vertices of the digraph) with $|A_i|=d^+_i$ for each $i$, and similarly $B$ partitioned into nonempty sets (bins) $B_i$, $i=1, \ldots, n$  with $|B_i|=d^-_i$ for each $i$. We write $\alpha(a_i)=j$ if $a_i\in A_j$, and $\beta(b_i)=j$ if $b_i\in B_j$. A random element of $\Psd$ is a random bijection $\phi:A\to B$ together with a random permutation $\sigma$ of $[n]$, such that the pair $(\phi,\sigma)$ is chosen u.a.r. Each element in $\Psd$ can be mapped in a natural way to a pairing in $\Pnm$, obtained by identifying points in $A_{\sigma(j)}$ and points in $B_j$ with out-points and in-points of bin $j$. This corresponds in turn to a multidigraph $M$ which has an arc $(u,v)$ for each point $a_i\in A_{\sigma(u)}$ such that $\phi(a_i)\in B_v$, or equivalently the arc (multi)set is $\{\sigma^{-1}(\alpha(a_i))\beta(\phi(a_i))\,:\, 1\le i
\le m\}$. Observe that $M$ has indegree sequence $\dmv$ and an outdegree sequence which is a random permutation of $\dpv$. As usual, all graph theory statements referred to an element in $\Psd$ should be understood in terms of the corresponding multidigraph.

Define $U$  to be the event, defined on any relevant
probability spaces, that there is a complex proper sink-set
containing at most $m/2$ arcs. Ultimately, we will do the calculations in the space $\Psd$ with $\dv$ randomised according to its distribution in the space $\hD$. Call this space $\Psnm$. Averaging over $\dv$ makes computations a little easier than arguing about its typical values. In fact, observe that the distribution of a random degree sequence $\dv\in\hD$ stays invariant if we randomly permute the entries of the outdegree sequence $\dpv$. Hence, we deduce that
\begin{equation}\lab{eq:PstoP}
\pr_\Pnm(U) = \ex_\hD \left(\pr_\Pd(U)\right) = \ex_\hD \left(\pr_\Psd(U)\right) = \pr_\Psnm(U).
\end{equation}
Thus, in view of Corollary~\ref{c:PtoG} and Lemma~\ref{l:useful}, we have
\begin{equation}\label{eq:PstoG}
\pr_\dico(U) = \pr_\Pnm(U\mid S) \le (1+o(1))e^{\lambda^2/2} \, \pr_\Pnm(U) = O(\pr_\Psnm(U)).
\end{equation}
Therefore, we only need to show that $\pr_\Psnm(U)=o(1)$ in order to prove the theorem statement for complex sink-sets. The result extends immediately to complex source-sets by considering the converse digraph.

The remainder of the proof consists of bounding the probability that an element of $\Psnm$ has a complex sink-set with at most $m/2$ arcs. Observe that, if $S$ is a complex sink-set and $v_0$ is a vertex in $S$ with outdegree strictly greater than $1$ (there must exist at least one of these because $S$ is complex), then the set $S'\subseteq S$ of vertices reachable from $v_0$ is also a complex sink-set.
Therefore, we only need to consider complex sink-sets which are precisely the set of vertices reachable from some vertex $v_0$.

Given a vertex $v_0$, the following algorithm will terminate with $S$ being the set of vertices reachable from $v_0$. The algorithm works by maintaining a set  $S$   of bins $A_i$ corresponding to 
vertices reachable from $v_0$, and investigating the vertices reachable from $S$. It does this by looking at the points in bin in $S$. The set $T$  contains precisely such points  which have not yet been investigated.

\bigskip

\noindent {\bf Algorithm}

\smallskip

\noindent
Let $v_0$ be the initially chosen vertex. Start with $S=\{v_0\}$,
$T=A_{\sigma(v_0)}$,
  and repeat the following until
$T$ is empty. Pick $a_i\in T$, add to $S$ the vertex
$v=\beta(\phi(a_i))$ (if it is not already there), delete
$a_i$ from $T$ and, if $v$ was not already in $S$, add all elements in
$A_{\sigma(v)}$ to $T$.

\bigskip

If the algorithm terminates with $S$ being a complex sink-set
containing at most half of the arcs of $M$, we say that it terminates {\em
properly}, and otherwise improperly.   We complete the proof of the theorem by showing that the probability that
there exists a vertex $v_0$ such that the algorithm terminates properly, when
begun from $v_0$, is $o(1)$.

As is common in analysing algorithms like this, we will make use of the fact that, conditioning on any set of values of a uniformly random permutation, the remaining values are still uniformly at random. Thus,  the algorithm can be performed
simultaneously with the generation of the random  bijection $\phi$ and
permutation $\sigma$. At the start, $\phi$ and $\sigma$ entirely undetermined and we can choose
$\phi(a_i)$
at random from the unused points of $B$ at each step of the algorithm.
Similarly, we may choose $\sigma(v_0)$ initially, and then
$\sigma(v)$ at each step where the vertex $v$ was not already in
$S$, randomly from the indices $i$ of the unused bins $A_i$. 
Thus, we may initially choose u.a.r.\ a permutation $\phi_1,\ldots, \phi_m$ of $B$, and independently a permutation $\sigma_1,\ldots, \sigma_n$ of $[n]$ u.a.r., and use $\phi_1$ for the first value of $\phi$ called for in the algorithm, $\phi_2$ for the second, and so on, and similarly for $\sigma$.  Set $K_k=\{\phi_1,\ldots, \phi_{k}\}$ and $J_s=\{\sigma_1,\ldots, \sigma_{s}\}$. Since the $\phi_i$ and $\sigma_i$  are pre-chosen randomly, it follows  that, for given $k$ and $s$,
\bel{modelb}
\mbox{for given $k$ and $s$,  $J_s\subseteq[n]$ and $K_k\subseteq B$ are independent and u.a.r.} 
\end{equation}
 In particular, the joint distribution of $J_s$ and $K_k$ does not   not depend on the algorithm, which is the important feature that simplifies analysis.

Now define
\begin{equation}\lab{khat}
\hat k =\sum_{j\in J_s} |A_{j}|,
\end{equation}
and let $U_{v_0}$ denote the event that the algorithm terminates properly,
with $k$ and $s$ defined as above, in particular with $S$ being a complex sink-set with at most $m/2$ arcs.
  In the event  $U_{v_0}$, since the
termination condition implies that $T$ is empty, it follows that
\begin{equation}\lab{condb}
\hat k=k.
\end{equation}
Also define
\begin{equation}\lab{shat}
\hat s =|\{u\,:\, u=v_0 \mbox{ or } K_k\cap B_u\ne
\emptyset\}|.
\end{equation}
Note that at each step, since $\beta (\phi(a_i))$ is added to $S$, we have
$S=\{v_0\} \cup \{u\, :\, K_k\cap B_u\ne \emptyset\}$ and hence
\begin{equation}\lab{conda}
\hat s =s.
\end{equation}
Moreover, the fact that $S$ is complex is equivalent to the condition that there are more arcs chosen than vertices in $S$, and so $k>s$. Hence, an upper bound on $\pr(U_{v_0})$ is the probability that~\eqref{condb} and~\eqref{conda} hold for some $k$ and $s$ with $k\le m/2$ and $s<k$. 
% (For such $k$ and $s$ to occur at the termination of the algorithm, 
% it  would  also be required that these conditions are not met for some 
% smaller values of $k$ and/or $s$, but we may ignore this.)

% Conditional on the values of $k$ and $s$,~\eqref{modelb} shows
% that~\eqref{condb} is independent of $\vd$, and~\eqref{modela} shows
% that~\eqref{conda} is independent of $\vdel$.
% Note that, given $k$ and $s$, the distribution of $\hat s$ depends only on 
% $\dmv$ and $K_k$, and that of $\hat k$ depends only on $\dmv$ and $J_s$, 
% and moreover these two variables $\hat s$ and $\hat k$ are independent.
Denote the event that~\eqref{conda} holds,
given $k$ and $s$, with $\hat s$ generated according to~\eqref{shat} given~\eqref{modelb}, by $H^-_{k,s}$, and similarly the event that~\eqref{condb} holds, given
$k$ and $s$, with $\hat k$ generated according to~\eqref{khat}, by $H^+_{k,s}$. Also put
%
% \begin{equation}\lab{Hksdef}
$ H_{k,s}=H^+_{k,s}\wedge H^-_{k,s}.$
% \end{equation}
%
We will prove that
\begin{equation}\lab{sumHks}
\pr\bigg(\bigcup_{k\le m/2}\,\bigcup_{s<k} H_{k,s}\bigg) = o(n^{-1}).
\end{equation}
Then $\pr_\Psnm(U_{v_0})  = o(n^{-1})$,
and the result follows by taking the union bound $\pr_\Psnm(U) \le \sum_{v_0}\pr_\Psnm(U_{v_0})$ and from~\eqref{eq:PstoG}.

It only remains to show~\eqref{sumHks}, for which we split $k$ into two
intervals.
\smallskip

\no {\em Case 1. $\log^4 n < k \le m/2$.}

\no We first bound probabilities in the distribution of
$\hat s$ as determined by~\eqref{shat}. Recall the truncated Poisson distribution as defined in Section~\ref{s:basics}.
Let $\Omega=\Omega(n,c,q)$ denote the probability space in
which there are $n$ bins $B_i$ with $|B_i|=d^-_i$, where $d^-_1, \ldots , d^-_n$ are independent random variables each with the distribution of $\tpol$, and such that a random subset $T$ of the  points in the bins is chosen by including each point independently with probability
\[
q = k/m.
\]
Let $\hat s$ be the number of the bins that are either occupied by at least one point of $T$ or happen to be the bin $v_0$. It follows from Lemma~\ref{l:binsfromballs} that
\begin{equation}
\lab{eq:snq}
\pr_\Omega\big(|\hat s-nq'|>\sqrt{nq'}\log n\big)=o(n^{-3}),
\end{equation}
where $q'$ is the probability that a bin contains some point of $T$. Note that $q'>q(1+\epsilon)$ for some positive constant $\epsilon$ that can be determined from~\eqref{eq:qprime}.
Now define $E$ to be the event in $\Omega$ that the total content of bins is $\sum_{i=1}^nd^-_i=m$ and that exactly $k$ points are chosen in $T$. Observe that, in the probability space $\Omega$ conditional upon $E$, the number $\hat s$ of bins containing at least one element of $T$ is distributed as in the definition of $H^-_{k,s}$. From Lemma~\ref{l:sum}, $\sum_{i=1}^nd^-_i=m$ holds with probability $\Theta(n^{-1/2})$, and---conditional on that---the event $|T|=k$ has probability $\Theta(k^{-1/2})$, since $|T|\distrib\bindis(m,k/m)$. Hence, $\pr_\Omega(E)=\Omega(n^{-1})$ and by~\eqref{eq:snq}
\begin{equation}\lab{eq:PHksm}
\pr\Bigg(\bigcup_{|s-nq'|>\sqrt{nq'}\log n}H^-_{k,s}\Bigg) = \pr_\Omega(|\hat s-nq'|>\sqrt{nq'}\log n\mid E) = o(n^{-2}).
\end{equation}
The next (and simpler) step is to define  $\Omega'=\Omega'(n,c,s)$ to be the probability space in which there are $n$ bins $A_i$ with $|A_i|=d^+_i$ all independent random variables each with the distribution of $\tpol$, and such that a uniformly random set of $s$ of the bins is chosen. Assume that $s$ lies in the range $|s-nq'|\le\sqrt{nq'}\log n$, and in particular $s=\Omega(\log^4 n)$.
Let $\hat k$ be the total number of points in the selected bins. From Lemma~\ref{l:concentration} we obtain the tail bound
\begin{equation}\lab{eq:kcs}
\pr_{\Omega'}\left(\hat k \le \frac{cs}{1+\epsilon/2}\right) \le e^{-\Theta(s^{1/3})} = o(n^{-4}).
\end{equation}
Now let $E'$ be the event in $\Omega'$ that $\sum_{i=1}^nd^+_i=m$, and recall from Lemma~\ref{l:sum} that $\pr_{\Omega'}(E')=\Theta(n^{-1/2})$. Observe that in $\Omega'$ conditional upon $E'$ the distribution of $\hat k$ is the same as the one in the definition of $H^+_{k,s}$. Moreover, the fact that $|s-nq'|\le\sqrt{nq'}\log n$ implies that $k\le cs/(1+\epsilon+o(1))$. In view of all that and from~\eqref{eq:kcs}, we obtain that for $|s-nq'|\le\sqrt{nq'}\log n$
\begin{equation}
\pr(H^+_{k,s}) \le \pr_{\Omega'}\left(\hat k \le \frac{cs}{1+\epsilon+o(1)}\mid E'\right) = o(n^{-3}).
\lab{eq:PHksp}
\end{equation}
Taking the union bound over all $s$ such that $|s-nq'|\le\sqrt{nq'}\log n$, combining it with~\eqref{eq:PHksm} and summing over all $k$ between $\log^4 n$ and $m/2$ completes the proof of~\eqref{sumHks} for this range of $k$.

\smallskip

\no {\em Case 2. $k\le\log^4 n$.}

\no For $\dv\in\hD$, the event $d^-_{\max}<\log^2n$, or equivalently $|B_i|<\log^2 n$ for each $i$, holds with probability $1-o(n^{-1})$. This follows readily from bounding the probability of the complement in $\cD$ and then conditioning upon $\sum_{i}d^-_i=m$ (see~\eqref{eq:tailbounds} and Lemma~\ref{l:sum}). Since $\Psnm$ is just $\Psd$ with $\dv$ distributed as in $\hD$, we may focus on $\Psd$ for a particular $\dv$ satisfying the above property. Note that $\hat s = s<k$ according to~\eqn{shat} if  the random choice $\{\phi_1,\ldots, \phi_{k}\}$ of elements of $B$ determines at most $k-2$ bins other than $v_0$. This has probability $O(k^4(\log^2 n/m)^2)$.
% this is thinking of choosing sequentially: for one collision choose two points and
% make the second go in the same bin as the first. Then square this.
% Another way: choose two bins (n^2) then prob of hitting one twice is k^2  log^4 n
% / m^2. THen square this. This is weaker.
Hence, in  $\Psd$
\[
\sum_{k\le \log^4 n}\pr\bigg(\bigcup_{s<k}H^-_{k,s}\bigg)= \sum_{k\le \log^4 n}\pr(\hat s \le k-1) = O(\log^{24}n/n^2)=o(1/n).
\qedhere
\]
\end{pf}

Next consider plain s-sets of $\dico$.
\begin{prop}\lab{p:plain} Suppose that $c=m/n$ is bounded and bounded away from $1$. The probability
that a digraph in $\dico$
has no plain s-set is asymptotic to
\begin{equation}\lab{nononcomplex}
\frac{e^\lambda(e^\lambda-1-\lambda)^2}
{(e^{2\lambda}-e^\lambda-\lambda)(e^\lambda-1)},
\end{equation}
with $\lambda$  determined by the equation
$c=\lambda e^\lambda/(e^\lambda-1)$.
\end{prop}
\begin{pf}
The simplest sink-sets or source-sets are those whose vertices induce a directed cycle. Call them \emph{sink-cycles} or \emph{source-cycles} accordingly. An \emph{s-cycle} is just a set of vertices which is either a sink-cycle or a source-cycle. Observe that each plain s-set must contain one s-cycle, so we can restrict our attention to s-cycles.

For any constant natural $k\ge1$, let $C_k$ be the number of s-cycles of order at most $k$. Let $D$ be the number of double arcs. Define
\[
\mu_k = \sum_{j=1}^k \frac{2(c/e^\lambda)^j-(c/e^{2\lambda})^j}{j}.
\]
Easy computations show that $2(c/e^\lambda)^j>(c/e^{2\lambda})^j$, so that there are no cancellations in any term of the definition of $\mu_k$.
We first claim that $\ex_\Pnm C_k \sim \mu_k$, $\ex_\Pnm D \sim \lambda^2/2$, and moreover $C_k$ and $D$ are asymptotically jointly independent Poisson. Elementary calculations show that
\[
\mu=\lim_{k\to\infty}\mu_k = \log\left(\frac{(e^{2\lambda}-e^\lambda-\lambda)(e^\lambda-1)} {e^\lambda(e^\lambda-1-\lambda)^2}\right).
\]
On the other hand, we claim that the probability there is an   s-cycle of order greater than $k$ can be bounded by some function $f_k$ such that $\lim_{k\to\infty}f_k=0$. In view of all this, setting $S$ to be the event that $\Pnm$ has no multiple arcs and $V$ to be the event in $\dico$ or $\Pnm$ that there are no s-cycles, we get
\[
\pr_\Pnm(V\cap S) \sim e^{-\mu-\lambda^2/2}.
\]
Then the proof of the result follows immediately from Lemma~\ref{l:useful}(b) and the fact that $\pr_\dico(V) = \pr_\Pnm(V\mid S)$.

Now we proceed to verify the claims we made about $C_k$, $D$ and the expected number of ``long" s-cycles. To make the computations easier, we generate the elements of $\Pnm$ using a slight variation of the $\Psnm$ model in which the in-points $a_1,\ldots,a_m$ (resp. out-points $b_1,\ldots,b_m$) are assigned independently and u.a.r.\ to the in-bins $A_1,\ldots,A_n$ (resp. out-bins $B_1,\ldots,B_n$) conditional upon each bin receiving at least one point (note that the degree sequence thus obtained is distributed as in $\hD$). In addition to that, a random bijection $\phi$ of the out- and in-points, and a random permutation of the labels of the out-bins are chosen as before independently and u.a.r.\ (alternatively we may consider $\sigma$ to be a random bijection of the out- and in-bins).

First we wish to compute the joint factorial moments of $C_k$ and $D$. We shall index all possible s-cycles of length at most $k$ by their position (i.e.\ the vertices they use in cyclic order). More precisely, the position of a cycle of length $\ell$ is determined by a tuple of $\ell$ distinct in-bins $B_{i_1},\ldots,B_{i_\ell}$ given in cyclic order together with and ordered tuple of out-bins $A_{i_1},\ldots,A_{i_\ell}$. A random element of $\Psnm$ has an s-cycle at $c$, if it has an s-cycle on vertices $v_1,\ldots,v_\ell$ where each vertex $v_j$ corresponds to the bins $A_{\sigma(i_j)}$ and $B_{i_j}$.

Fix $c_1,\ldots,c_r$, where each $c_i$ is the position of a cycle of length $\ell_i$ and the bins used for each position are pairwise disjoint. Let $X_{c_1,\ldots,c_r}$ be the indicator function for the event that there is an s-cycle at each position $c_i$. We compute the probability that this event holds. The probability that the out-bins are assigned to the corresponding in-bins is
\begin{equation}
1/[n]_{\ell_1+\cdots+\ell_r} \sim 1/n^{\ell_1+\cdots+\ell_r}.
\label{eq:Pscycles0}
\end{equation}
Condition on this, and note that the degrees of the bins and the matching of the points occur independently from that. Now we claim that the probability that the right s-cycles occur at $c_1,\ldots,c_r$ is asymptotic to
\begin{equation}
  \prod_i (2a^{\ell_i}/n^{\ell_i}-a^{2\ell_i}/m^{\ell_i}),
\label{eq:Pscycles}
\end{equation}
where $a=\lambda/(e^\lambda-1)=c/e^\lambda$. Observe that the events of having a sink-cycle or having a source-cycle at $c_i$ are not disjoint, so the probability of the union is the sum of probabilities minus the probability of having both a sink- and a source-cycle at $c_i$. Thus, in order to estimate~\eqref{eq:Pscycles}, we can specify for each $c_i$ one of the three former events (sink-cycle, source-cycle or both). More precisely, given $r_1+r_2+r_3=r$ and relabelling $c_1,\ldots,c_r$ to $c_1,\ldots,c_{r_1},c'_1,\ldots,c'_{r_2},c''_1,\ldots,c''_{r_3}$ and $\ell_1,\ldots,\ell_r$ to $\ell_1,\ldots,\ell_{r_1},\ell'_1,\ldots,\ell'_{r_2},\ell''_1,\ldots,\ell''_{r_3}$, we shall compute w.l.o.g.\ the following probability: having a sink-cycle at $c_1,\ldots,c_{r_1}$ (and possibly a source-cycle too); having a source-cycle at $c'_1,\ldots,c'_{r_2}$ (and possibly a source-cycle too); or having both a sink- and a source-cycle at $c''_1,\ldots,c''_{r_3}$. We shall see that this probability is asymptotic to
\begin{equation}
  \frac{a^{\ell_1+\cdots+\ell_{r_1}+\ell'_1+\cdots+\ell'_{r_2}+2(\ell''_1+\cdots+\ell''_{r_3})}}{n^{\ell_1+\cdots+\ell_{r_1}+\ell'_1+\cdots+\ell'_{r_2}} m^{\ell''_1+\cdots+\ell''_{r_3}}},
\label{eq:Pscycles2}
\end{equation}
and this leads to~\eqref{eq:Pscycles} by easy inclusion-exclusion. 

To arrive at~\eqref{eq:Pscycles2}, we require that $\ell_1+\cdots+\ell_{r_1}+\ell''_1+\cdots+\ell''_{r_3}$ specific out-bins determined by $c_1+\cdots+c_{r_1}+c''_1+\cdots+c''_{r_3}$ contain exactly one point each. By symmetry, the probability of this is $(\ex [N]_{\ell_1+\cdots+\ell_{r_1}+\ell''_1+\cdots+\ell''_{r_3}})/[n]_{\ell_1+\cdots+\ell_{r_1}+\ell''_1+\cdots+\ell''_{r_3}}$, where $N$ is the number of out-bins with only one point. $N$ is concentrated around $an$ by~\cite[Lemma~1]{CW} since the distribution of balls in bins is truncated multinomial.
%query: note this reference used here 
Hence $\ex [N]_{\ell_1+\cdots+\ell_{r_1}+\ell''_1+\cdots+\ell''_{r_3}} \sim (an)^{\ell_1+\cdots+\ell_{r_1}+\ell''_1+\cdots+\ell''_{r_3}}$ and the probability is $(1+o(1)) a^{\ell_1+\cdots+\ell_{r_1}+\ell''_1+\cdots+\ell''_{r_3}}$. Analogously, we need that some specific $\ell'_2+\cdots+\ell'_{r_2}+\ell''_1+\cdots+\ell''_{r_3}$ in-bins contain only one point, which is independent from the previous and has probability $(1+o(1)) a^{\ell'_2+\cdots+\ell'_{r_2}+\ell''_1+\cdots+\ell''_{r_3}}$. 
Conditional upon all this, we need that for each $c_1,\ldots,c_{r_1}$, the only point in each out-bin is matched to some point in the corresponding in-bin; for each $c'_1,\ldots,c'_{r_2}$, the only point in each in-bin is matched to some point in the corresponding out-bin; and for each $c''_1,\ldots,c''_{r_3}$, the only point in each out-bin is matched to the only point in the corresponding in-bin. 
Observe that the number of points in these out-bins that have not been exposed remains independent truncated Poisson conditional to fixed sum $m-\ell_1+\cdots+\ell_{r_1}+\ell''_1+\cdots+\ell''_{r_3}$. An analogous thing happens for in-bins that were not exposed and the sum of their degrees is $m-\ell'_1+\cdots+\ell'_{r_2}+\ell''_1+\cdots+\ell''_{r_3}$. The probability of matching the in- and out-points appropriately for s-cycles at $c''_1,\ldots,c''_{r_3}$ is $1/[m]_{\ell''_1+\cdots+\ell''_{r_3}} \sim 1/m^{\ell''_1+\cdots+\ell''_{r_3}}$. 
We condition on that and on the event that no out-point corresponding to $c_1,\ldots,c_{r_1}$ is matched to any in-point corresponding to $c'_1,\ldots,c'_{r_2}$ (which happens with probability $1+o(1)$). This makes the construction of the remaining sink-cycles independent from that of source-cycles. 
For the sink-cycles, we have to match the only point in each out-bin with some point in the corresponding in-bin. By symmetry, the first matching has probability $1/(n-\ell'_1+\cdots+\ell'_{r_2}+\ell''_1+\cdots+\ell''_{r_3})\sim 1/n$. Conditional to some matchings being exposed, the probability that the next out-point is matched to a point in an in-bin which contains one matched point already is $O(1/n)$ since there is negative correlation between the events that two out-points are matched to in-points in the same bin (condition on any given degree). 
Thus that out-point is matched to some point in an unexposed in-bin with probability $1+o(1)$ and conditional to that, again by symmetry, chooses the right in-bin with probability $(1+o(1))/n$. This gives a probability   $(1+o(1))/n^{\ell_1+\cdots+\ell_{r_1}}$ for having the matchings required for the sink-cycles. Analogously, the source-cycles give a $(1+o(1))/n^{\ell'_1+\cdots+\ell'_{r_2}}$ factor, and this establishes the estimate~\eqref{eq:Pscycles2}.

So far we were dealing with fixed cycle tuples $c_1,\ldots,c_r$. Let $C_k$ be the random number of s-cycles occurring. To compute the $r$-th factorial moment it suffices to multiply~\eqref{eq:Pscycles0} and ~\eqref{eq:Pscycles} by the number of ways of choosing $r$ different $c_1,\ldots,c_r$, which is
\[
\sum_{\ell_1,\ldots,\ell_r\in\{1,\ldots,k\}}\frac{([n]_{\ell_1+\cdots+\ell_r})^2}{\ell_1\cdots\ell_r}.
\]
Hence,
\[
\ex[C_k]_r \sim \left( \sum_{\ell=1}^k \frac{2a^{\ell_i} - (a^2/c)^{\ell_i}}{\ell_i} \right)^r = {\mu_k}^r.
\]

Let $D$ be the number of double arcs occurring. Recall that the out-points are placed in the out-bins (and the in-points in the in-bins) uniformly at random and independently conditional upon getting at least one point in each bin. We index double arcs according to their position, where each position $j$ is a set of two different out-points in the same out-bin along with a set of two different in-points in the same in-bin. Let $Z$ be the number of positions for double arcs. We have the trivial bound $Z\le m^4$. Combining together Lemmas~\ref{l:multinomial}, \ref{l:variance} and~\ref{l:sum}, we have that $|Z-(\lambda c n/2)^2|<n^{1.6}$ with probability $1-O(n^{1/2}e^{-\log^3n})$. Hence, $\ex[Z]_s\sim(\lambda c n/2)^{2s}$. We say that there is a double arc at $j$ if the out-points are matched to the in-points in any of the two possible ways. Fix $s$ different positions $j_1,\ldots,j_s$. The probability of having double arcs at $j_1,\ldots,j_s$ is is $2^s/[m]_{2s}$. Therefore,
\[
\ex[D]_s = \ex[Z]_m 2^s/[m]_{2s} \sim (\lambda^2/2)^s.
\]
In order to compute joint moments of $C_k$ and  $D$, we condition to s-cycles happening at some fixed positions $c_1,\ldots,c_r$ and we specify the type of each cycle (sink-, source- or both) in the same fashion we used in the previous computations (sink does not exclude source and vice-versa). To make computations easier we also condition on the particular in-points and out-points matched to create the s-cycles. Conditional upon all this, we compute $\ex[D]_s$. The same computations we did before are still valid if applied to the in- and out-points that were not used in the construction of the s-cycles, and this yields the the same asymptotic value $(\lambda^2/2)^s$. So
\[
\ex[C_k]_r[D]_s \sim \left( \sum_{\ell=1}^k \frac{2a^{\ell_i} - (a^2/c)^{\ell_i}}{\ell_i} \right)^r (\lambda^2/2)^s
\]
The claim that the distributions are asymptotically independent   Poisson now follows by the standard method of moments.

It only remains to bound the probability of existence of s-cycles of length greater than $k$ by some function $f_k$ such that $\lim_{k\to\infty} f_k=0$. It is enough to deal with sink-cycles, since the result for source-cycles follows by considering the converse digraph. Take a length $\ell>k$. We now condition on the number  $N$ as defined above.
 We can choose $([n]_\ell)^2/\ell$ different positions for such a cycle. For each of these, the probability that the bins are matched the right way is $1/[n]_\ell$ (regardless of $N$). The probability that each of the $\ell$ out-bins contains exactly one point is at most $ (N/n)^\ell$. The probability that each of the $\ell$ out-points is matched to a point in the corresponding in-bin is at most $1/[n]_\ell$ (since conditional upon $i$ matched pairs, the probability of the next matched pair is $1/(n-i)$ times the probability of hitting a point in an in-bin not previously hit). This is again regardless of $N$.
 
Putting this together, this expectation is at most
\[
\sum_{\ell>k}  (N/n)^\ell /\ell.
\]
this tends to $0$ for large $k$ provided $N/n<(a+1)/2$. The probability that $N$ is larger than this is $o(1)$ by the concentration mentioned above.
%%%%
% Note - that lemma applies because m/n is bounded away from 1.
%%%%
\end{pf}
{From} Theorem~\ref{t:dicoreenum}, Proposition~\ref{p:complex} and Proposition~\ref{p:plain} we immediately obtain Theorem~\ref{t:main} for the case $c$ is bounded and bounded away from $1$.
%
%
%-----------------------------------------
%-----------------------------------------
%
\section{Very sparse case: $c\to 1$}\lab{s:vsscd}
Here we proof Theorem~\ref{t:main} for the case $c\to1$. Thus $r=m-n=o(n)$, and we assume $r=m-n\to\infty$.  
%
%
%\subsection{The heart model}\lab{s:hearts}
%
We define the directed graph analogue of the kernel of a graph as follows.  A cycle component of a digraph is a connected component which is simply a directed cycle. 
A digraph with each in-and outdegree at least $1$ and with no cycle components is called a \emph{preheart}.
The {\em heart} of  a preheart $G$ is the multidigraph $H(G)$
obtained from $G$ by repeatedly choosing a vertex
$v$ of in- and outdegrees both $1$, deleting $v$ and its two incident arcs $uv$ and $vw$, and inserting the arc $uw$. The condition that
$G$ contains no isolated cycle ensures that the heart is always a multidigraph.  The vertices of $H(G)$ are just the vertices of $G$ of total degree at least 3.

Note that a digraph is strongly connected iff it is an isolated cycle or a preheart with strongly connected heart. Thus we may use ideas similar to those in Section~\ref{s:sscd} to study the heart, as a key step to enumerate strongly connected digraphs.
\remove{ %%%%%%%%%%%%%%%%%%%%%%%%%
Connectivity properties of the heart are perhaps no easier to prove than for just the $(1,1)$-dicore. Consider for instance the probability that a heart has no complex s-set. Since many vertices may have outdegree $1$, it is hard to see that this is any easier than for the $(1,1)$-core unless attention is paid to the matching of in- and outdegrees. We will take care of this by combining the idea of heart and a randomisation of the in- and outdegree sequences, similarly as in the $\Pd$ and $\Psd$ models in Section~\ref{s:sscd}.
}%%%%%%%%%%%% end remove
Connectivity properties of the heart 
% CHECK!!!! changed meaning to opposite sense!!
can be easier to prove than for just the $(1,1)$-dicore. In the dicore, some complex s-sets can involve  many vertices of in- and outdegree 1 and just a few other vertices. We will focus on the heart and also use randomisation of the in- and outdegree sequences,   as in the $\Pd$ and $\Psd$ models in Section~\ref{s:sscd}.

Consider any given degree sequence $\dv\in\hD$, and let $T=T(\dv)=\{i:d^+_i+d^-_i\ge 3\}$. We put $n'=|T|$ and
$m'=\sum_{i\in T} d^+_i=\sum_{i\in T} d^-_i$, and note that $m-n=m'-n'$.  For simplicity of presentation, renumber the vertices if necessary so that $T=[n']$.

Let $\Hd$ be the probability space of \emph{heart configurations} generated as follows. For each $i\in T$ consider consider a bin containing labelled points of two types, namely $d^+_i$ \emph{out-points} and $d^-_i$ \emph{in-points}, and then choose a random matching of the in-points with the out-points (there are $m'$ of each kind). Note that each heart configuration in $\Hd$ corresponds to a multidigraph on vertex set $T$ obtained in a natural way by identifying bins with vertices and adding an arc $(u,v)$ for each out-point in $u$ matched to an in-point in $v$.

Moreover, given a heart configuration $H$, we construct a \emph{preheart configuration} $Q$ by taking an assignment of $[n]\setminus T$ to the arcs of $H$ (i.e. the pairs of matched up points), such that the numbers assigned to each arc are are given a linear ordering. Denote this assignment, including the linear orderings, by $f$. Let $\Qd$ be the probability space of random preheart configurations created by taking $H\in\Hd$ and choosing $f$ u.a.r. Note that each $Q\in\Qd$ corresponds to a multidigraph with $n$ vertices, $m$ arcs and degree sequence $\dv$. Henceforth, any graph terminology referring to a heart or preheart configuration should be interpreted in terms the corresponding multidigraph.
\begin{lem}\lab{l:heartuniform}
The digraphs generated from the restriction of $\Qd$ to simple preheart configurations (i.e.\ with no multiple arcs) are uniformly distributed.
\end{lem}
\begin{pf}
Each simple digraph comes from $\prod_{i=1}^n d^+_i!d^-_i!$ different preheart configurations.
\end{pf}
As will become apparent later in the argument, it turns out that the degree sequence distribution induced by the uniform probability space of all prehearts on $n$ vertices and $m$ arcs is close in some sense to  $\hD$. This motivates considering the probability spaces $\Hnm$ and $\Qnm$, defined by choosing a random element from $\Hd$ and $\Qd$ respectively, where the degree sequence $\dv$ is also random and distributed as in $\hD$.

Given a degree sequence $\dv\in\hD$, we distinguish four kinds of vertices depending on whether their in- and outdegree are equal to 1, or larger. For $i,j\in\{1,2\}$, let $N_{i,j}$ be the set of vertices with indegree of type $i$ and outdegree of type $j$ (type $1$ means $1$ and type $2$ means greater than 1). Let $\ba=(a_{1,1},a_{1,2},a_{2,1},a_{2,2})$, where $a_{i,j}=|N_{i,j}|$. This $\ba$ is of course a function of $\dv$.
Observe that any $\ba$ which is \emph{feasible} (i.e.\ occurs in $\hD$ with nonzero probability) satisfies $a_{1,1}+a_{1,2}+a_{2,1}+a_{2,2}=n$, $1\le a_{1,2}+a_{2,2}\le r$ and $1\le a_{2,1}+a_{2,2}\le r$. Conversely, it is easy to check that, for sufficiently large $n$, any nonnegative tuple $\ba$ satisfying the above conditions is feasible.
Note also that   $n'=a_{1,2}+a_{2,1}+a_{2,2}$ and $m'=r+a_{1,2}+a_{2,1}+a_{2,2}$.

We will want to condition on ``typical" values of $\ba$. Denote by $\good$ the event that
\[
|a_{1,2}-r| \le \sqrt r\log r, \quad
|a_{2,1}-r| \le \sqrt r\log r \quad\text{and}\quad
a_{2,2} \le \max\{2r^2/n,\sqrt r\}.
\]
Note in particular that $\good$ implies 
\bel{goodcond}
n'\sim 2r\qquad m'\sim 3n'/2\sim 3r,\qquad a_{1,2}\sim a_{2,1}\sim r,\qquad a_{2,2}=o(r).
\ee
 We next show something somewhat stronger than $\pr_\hD(\good)=1-o(1)$.
\begin{lem}\lab{l:gooda}
\[
\ex_\hD(m'(1-1_\good))=o(1).
\]
\end{lem}
\begin{pf}
First we observe that $m'$ is deterministically at most $3r$ in $\hD$. This upper bound is immediate from the fact that the underlying undirected graph of the heart has $n'$ vertices, $m'=n'+r$ edges and average degree $2m'/n'\ge3$.
Hence, by Lemma~\ref{l:useful}(a), it suffices to bound the probability that $\good$ fails by $o(1/r^2)$ in  $\cD$. Here,   $a_{1,2}$, $a_{2,1}$ and $a_{2,2}$ are binomially distributed with  expectations $r$, $r$ and $r^2/n$ respectively. (Note that $r\to\infty$, but $r^2/n$ need not be large.) 
Hence, standard bounds (if $r$ grows very slowly,~\eqn{chernsmallp} does not suffice, but in any case we can simply consider ratios of consecutive binomial probabilities) shows that the conditions on $a_{1,2}$ and $a_{2,1}$ in the definition of $\Gamma$ hold with probability $1-o(1/r^2)$. A similar argument ensures that $a_{2,2}$ has the required concentration with probability $1-o(1/r^2)$, but the analysis is split into two cases. If $r\le n^{3/5}$, then $r^2/n\le r^{1/3}$ and  we easily bound the probability that $a_{2,2}>\sqrt r$, for instance by comparing with a binomial with mean $r^{1/3}$. On the other hand, if $r>n^{3/5}$, we bound the probability that $a_{2,2}>2r^2/n$ using~\eqn{chernsmallp}.
\end{pf}
This result allows us to condition on feasible $\ba$ satisfying $\good$. In fact, for any given feasible tuple $\ba$, we denote by $\Ha$ and $\Qa$ respectively the probability spaces $\Hnm$ and $\Qnm$ conditional on having that particular $\ba$.
\begin{lem}\lab{l:heartnocomplex}
Let $\ba$ be any feasible tuple satisfying $\good$. Then a random heart configuration in $\Ha$ a.a.s.\ has no complex s-set of at most $m'/2$ arcs.
\end{lem}
\begin{pf}
 The argument shares many features with the proof of Proposition~\ref{p:complex}, in particular using auxiliary randomisation to simplify computations.  Let $N'$ denote $[n']$ (which was also $T$, the relevant set of vertices for the heart configuration).  For each $\dv\in\hD$ consider, as in the definition of $\Psd$, two sets of points $A=\{a_1, \ldots, a_{m'}\}$ and $B=\{b_1, \ldots , b_{m'}\}$, partitioned respectively into bins $A_1,\ldots,A_{n'}$ and bins $B_1,\ldots,B_{n'}$, with $|A_i|=d^+_i$ and $|B_i|=d^-_i$ for each $i\in N'$. We write $\alpha(a_i)=j$ if $a_i\in A_j$, and $\beta(b_i)=j$ if $b_i\in B_j$. Define the probability space $\Hsd$ to be a random bijection $\phi:A\to B$ chosen u.a.r.\ together with two random permutations $\sigma$ and $\tau$ of $[n']$, chosen independently of $\phi$ and of each other and u.a.r.\ subject to the conditions that $d^+_{\sigma(i)}=1$ whenever $d^+_i=1$, and $d^-_{\tau(i)}=1$ whenever $d^-_i=1$. We need an appropriate randomisation of the degrees. Thus, consider the probability space $\Hsa$, whose elements are selected at random from $\Hsd$ with $\dv$ a random member of $\hD$ but conditional on 
\remove{%%%%%%%%
the particular values of $\ba$.
}%%%%%%%%%
the particular value of the vector $\ba=\ba(\dv)$.

Observe that each element $H'$ in $\Hsa$ corresponds in a natural way to an element $H$ in $\Ha$, obtained by identifying the points in $A_{\sigma(j)}$ and those in $B_{\tau(j)}$ with the out-points and in-points, respectively, of bin (vertex) $j$ (in the same way that elements in $\Hsd$ can be mapped to elements in $\Hd$). Moreover, the $H$ obtained this way has the same distribution as in $\Ha$, since the distribution of the degree sequence and thus $\ba$ stay invariant after permuting the indices of the vertices in $N'$ by $\sigma$ and $\tau$ (so it does not matter if we condition to a particular $\ba$ before or after applying $\sigma$ and $\tau$).
% query : In next revision   should we define it the other way round?
 Hence, setting $U$ to be the event in $\Ha$ or $\Hsa$ that there is a complex sink-set containing at most $m'/2$ arcs, we have
\[
\pr_\Ha(U)=\pr_\Hsa(U).
\]
Henceforth we can do all calculations in $\Hsa$, which simplifies the analysis as   $\Psnm$ did in Section~\ref{s:sscd}.

By the same argument as in the proof of Proposition~\ref{p:complex}, in order to bound the probability of $U$, we can restrict our attention to complex sink-sets whose vertices are all reachable from some vertex $v_0$. If the set of vertices reachable from vertex $v_0$ is a complex sink-set then essentially the same algorithm as in
Section~\ref{s:sscd} will terminate with $S$ being such a sink-set.
We restate the algorithm in the current setting as follows:
\bigskip

\no
Start with $S=\{v_0\}$, $R=A_{\sigma(v_0)}$,  and repeat the following until $R$ is empty. Pick $i\in R$, add to $S$ the vertex $v$ such that $\phi(a_i)\in B_{\tau(v)}$ (if it is not already there), delete $i$ from $R$ and, if $v$ was not already in $S$, add all elements in $A_{\sigma(v)}$ to $R$.
\bigskip

As in the proof of Proposition~\ref{p:complex}, the algorithm
can be performed simultaneously with the generation of the random  bijection
$\phi$ and permutations $\sigma$ and $\tau$, piecemeal at each step of the
algorithm.

We need some notation to describe the generation of $\sigma$ and
$\tau$.  Let
\[
N_2^+=N_{2,1}\cup N_{2,2}=\{i\in N'\,:\,d^+_i>1\}, \quad N_2^-=N_{1,2}\cup N_{2,2}=\{i\in N'\,:\,d^-_i>1\},
\]
\[
N_1^+=N'\setminus N_2^+, \quad N_1^-=N'\setminus N_2^-.
\]
Also, at the start generate u.a.r.\ random permutations $\hat \sigma_j$ of $N_j^+$  and $\hat \tau_j$ of $N_j^-$ ($j=1,2$), and $\hat \phi$ of $B$, which we will view precisely as orderings of these sets (as in the proof of Proposition~\ref{p:complex}). Initially, let $\sigma(v_0)$ be the first element of $\hat\tau_j$, where $j$ is determined by $v_0\in N_j^+$. At each step, $\phi(a_i)$ is defined to be the next element of $B$ in the ordering $\hat \phi$. At each step where $\tau^{-1}(\beta(\phi(a_i)))$ has not yet been determined, choose $v$ to be the next unused member of  $N_j^-$ in the ordering $\hat \tau_j$. Set $\tau(v)=\beta(\phi(a_i))$. Then, if $\sigma(v)$ is not yet determined, define it as the next member of $N_j^+$ (where $v\in N_j^+$ determines $j$) in the ordering $\hat \sigma_j$.

%At any given stage, let $K_1\subseteq B$ denote the points $\phi(i)$ which
%have been chosen so far and for which
%$\beta(\phi(a_i))\notin N_2^-$, and  $K_2$ those for which
%$\beta(\phi(a_i))\in N_2^-$, and put $k_i=|K_i|$, $i=1$ and $2$.
At any given stage, when $k$ points $\phi(a_i)$   have been chosen so far,  let $K\subseteq B$ denote the set of these points, which must be the first $k$ points of $\hat \phi$. (This corresponds to the set $K_k$ in the proof of Proposition~\ref{p:complex}; we suppress the indices such as $k$ for simplicity.)  Also let $J^+$ denote the set of values $\sigma(v)$ determined so far (note
this is precisely $\{\sigma(v)\,:\, v\in S\}$), and somewhat asymmetrically, define $J^-$ to be the set of vertices whose image under $\tau$ has been determined. Then $J^-=S$ if $v_0$
was chosen at some stage as $v$, and otherwise $J^-=S\setminus v_0$.
Define the following random sets referring to a step after which   precisely
$k<m'$ arcs have been exposed:
\bean
J_1^+ &=&J^+\cap N_1^+, \\
J_2^+&=&J^+\cap  N_2^+,\\
J_1^- &=&J^-\cap N_1^-, \\
J_2^-&=&J^-\cap  N_2^-,
\eean
and put $t_1^+=|J_1^+|$ and so on.
%% It's a bit ugly at this point that $J^+$ is a set of half-bins (in image of sigma)
%%  but $J^-$ is a set of vertices (in preimage of tau).
Then at each step of the algorithm, conditional upon having given cardinalities that can feasibly occur, the permutations $\hat \sigma$ etc.\ determine these sets, and ensure that
each of these sets occurs u.a.r.\ as subsets of $N_1^+$,
$N_2^+$, $N_1^-$ and $N_2^-$ respectively, and the same property
holds for $K$ as a subset of points in $B$ with cardinality
$k$. Furthermore, all these sets occur jointly independently of each other. For $\tb = (t_1^+,t_2^+, t_1^-, t_2^-)$, let $\Omega(k,\tb  )$ denote the probability space of such independently chosen sets, $K$ and the $J_i^+$ etc., with these cardinalities.
Next define 
%(following the proof of Proposition~\ref{p:complex})
\bea
\hat k &=& \sum_{j\in J^+} |A_{j}|,  \lab{khat2}\\
\hat t_1^+&=&|(J^-\cup\{v_0\})\cap N_1^+|, \lab{hatt1p}\\
\hat t_2^+&=&|(J^-\cup\{v_0\})\cap  N_2^+|, \lab{hatt2p}\\
\hat t_1^-&=&|i\in N_1^-\,:\,K\cap B_i\ne\emptyset|, \lab{hatt1m}\\
\hat t_2^-&=&|i\in N_2^-\,:\,K\cap B_i\ne\emptyset|.
\lab{hatt2m}
\eea
By the form of the algorithm, at each iteration, precisely after the point when a new image of $\sigma$ is exposed, we have that $t_1^++t_2^+=t_1^-+t_2^-$ and also
\bel{newa}
\hat t_i^+=t_i^+ \ \mbox{and} \ \hat t_i^-=t_i^-,\quad i=1, 2.
\ee
Moreover, in the event $U_{v_0}$ that the algorithm terminates
with $S$ being a sink-set, we have
\bel{newb}
\hat k =k
\ee
and
\bel{newcondc}
k>t_1^++t_2^+
\ee
if it is complex.

\remove{ %%%%%%%%%%%%% remove

Then, setting $H_{k,\tb}$ to be the event that~\eqref{newa} and~\eqref{newb} hold, with $\hat k$, $\hat t_i^+$ and $\hat t_i^-$ defined by~\eqref{khat2}--\eqref{hatt2m} where $J^+$, $J^-$ and $K$ are chosen u.a.r.\ from $N'$, $N'$ and $B$ subject to $|J^+|=t_1^++t_2^+$, $|J^-|=t_1^-+t_2^-$ and $|K|=k$, we can easily 
find the coupling
\bel{newHv0}
  \pr_\Hsa(U_{v_0}) \le \pr\bigg(\bigcup_{k\le m/2}\ \bigcup_{t_1^++t_2^+=t_1^-+t_2^-<k}
H_{k,\tb}\bigg).
\ee
} %%%%%%%%%%end remove

Thus, setting $F_{k,\tb}$ to be the event that the tuple $\tb$ occurs in the algorithm after $k$ arcs are exposed, and $H$ the event that~\eqref{newa} and~\eqref{newb} hold, we have
by the union bound
\bel{event}
 \pr_\Hsa(U_{v_0})\subseteq \sum_{k\le m'/2}\ \sum_{t_1^++t_2^+=t_1^-+t_2^-<k}
 \pr_\Hsa(F_{k,\tb}\cap H ).
\ee
Note that   $H $ in also defined the space $\Omega(k,\tb  )$.

We now note  that, using the earlier observation that motivated defining $ \Omega(k,\tb  )$,    
$$ 
\pr_\Hsa(F_{k,\tb}\cap H )\le  \pr_\Hsa( H \mid F_{k,\tb})
=\pr_{\Omega(k,\tb  )}( H  ). 
$$
Thus it suffices to show 
\bel{newsumHks}
\sum_{k\le m'/2}\,\sum_{t_1^++t_2^+=t_1^-+t_2^-<k }\pr_{\Omega(k,\tb  )}( H  ) =
o(1/n'),
\ee
as the lemma follows from this using the argument in Proposition~\ref{p:complex} from~\eqref{sumHks} onwards.
 
Conditional on the values of $k$ and $\tb$, the random variables $\hat k$ etc.\ depend only on the random permutations $\hat \phi$ etc., and in particular the distribution of $\hat k$ only depends on $t_1^+$, $t_2^+$ and $\dpv$; the distributions of $\hat t_1^-$ and $\hat t_2^-$ only depend on $k$ and $\dmv$; the distributions of $\hat t_1^+$ and $\hat t_2^+$ only depend on $t_1^-$, $t_2^-$ and $\dpv$.

\noindent
{\em Case 1.} $ \log^4 n'<k\le m'/2$ 

Let $g=1/1000$. Let $E_1$ be the event that $|\hat t_1^- /(k/3)-1|\le g$ and $\hat t_2^-/(k/2)-1>-g$.
% query: subtracted g here. It seems if k = m' then we expect k/2 bins in N_2.
 Let $E_2$ be the event that $|\hat t_1^+ / t_2^- -1|\le g$ and $|\hat t_2^+ / t_1^- -1|\le g$. Let $E_3$ be the event that $|\hat k / (t_1^+ + 2t_2^+) -1|\le g$. Given any  fixed values for $k$ and $\tb$  with $k>\log^4 n'$, if both~\eqref{newa} and~\eqref{newb} hold then, clearly, at least one of  $E_1$, $E_2$  and $E_3$ must fail for $n'$ sufficiently large.  
\remove{ %%%%%%%%%%%%%%%%%
\begin{align*}
\hat t_1^+ &\sim t_2^-&
\hat k &\sim t_1^+ + 2t_2^+&
\hat t_1^- &\sim k/3\\
\hat t_2^+ &\sim t_1^-&&&
\hat t_2^- &> k/2\\
\end{align*}
} %%%%%%%%%%%%%%%%%%
We claim that each of $\overline E_1$, $E_1\setminus E_2$, and $(E_1\cap E_2)\setminus E_3$ have probability   $o((n')^{-5})$ in the spaces $\Omega(k,\tb  )$ occurring in~\eqref{newsumHks}. Thus $\pr_{\Omega(k,\tb  )}( H  ) =o((n')^{-5})$ in all cases,  yielding~\eqref{newsumHks} by summing over $k$ and the constrained $\tb$.

To verify the claims about the $E_i$, the same type of argument as in Case~1 in the proof of Proposition~\ref{p:complex} suffices. For instance, regarding $E_1$, recall that $K$ is a random subset of the points in $B$ of cardinality $k$. We can instead assume that the points of $B$ are independently choosen with probability $k/m'$, and condition later on obtaining precisely $k$ points, which holds with probability $\Theta(1/\sqrt k)$. Therefore, it is enough to show that $E_1$ has probability $1-o((n')^{-6})$ in the unconditional probability space where elements of $B$ are chosen independently. Noting  by~\eqref{goodcond} that $|N_1^-|\sim r$, $|N_2^-|\sim r$ and $m'=\Theta(n')$, we have by~\eqn{chernsmallp} that the probability that the number $t_1^-$ of points chosen in $N_1^-$ satisfies $|t_1^-/(k/3)-1|> g$ is $o(1/(n')^6)$. Similarly, from Lemma~\ref{l:binsfromballs} (applied to $|N_2^-|$ copies of $\tpo_2(\lambda)$ with $q=k/m'$), the probability that $| t_2^- /(k/2)-1|> g$ is $o(1/(n')^6)$. For $E_2$, note that $\hat t_1^+=|V\cap  N_1^+|$, where $V$ denotes the   set of the first  $t_2^-$ elements of $\hat \tau_2 $. Since $V$ is a random subset of 
$ N_2^-$, and since by~\eqref{goodcond} $| N_2^-\setminus N_1^+| = o(| N_2^-|)$, we have $|\hat t_1^+ / t_2^- -1|\le g$ with   probability $o((n')^{-5})$ provided say $t_2^- >\log^2 n'$. This is guaranteed by $E_1$.  The other statement in $E_2$ works exactly the same, and thus the probability of $E_1\setminus E_2$ is $o((n')^{-5})$.
Finally, for $E_3$, conditional on $\ba$, we just  consider the fixed number $r+a_{2,1}+a_{2,2}$ of   balls   thrown randomly into the $a_{2,1}+a_{2,2}$ bins conditional on at least two in each bin, (one ball in all other bins)   and argue as for~\eqref{eq:kcs} to deduce that  when $t_2^+$ bins are selected u.a.r., with high probability they contain approximately $2t_2^+$ balls.

\noindent
{\em Case 2.} $k\le \log^4 n'$ 

The argument for Case~2 in the proof of Proposition~\ref{p:complex} applies almost directly to the current setting, with of course $\Psd $ and $\Psnm$ replaced by $\Hsd$ and $\Hsa$. The only twist is that we have to show that, conditional upon $\ba$, the indegree  sequence  has maximum less than $\log^2 n$ with probability $1-o(1/n)$.  Such a sequence can be generated by putting $r+n'-a_{2,1}$ elements randomly into   $a_{1,2}+a_{2,2}$ bins subject to each bin receiving at least two balls. By~\eqn{goodcond} the excess of balls over bins is $o(r)$ and so the required property follows easily.% from [[[REFER TO, WITH SOME MODIFICATION, LEMMA USED AS IN CASE 2 OF EARLIER PROOF]].
\end{pf}
\begin{lem}\lab{l:heartsc}
Let $\ba$ be any feasible tuple satisfying $\good$. Then a random preheart configuration in $\Qa$ is simple and strongly connected with probability $1/9+o(1)$.
\end{lem}
\begin{pf}
A preheart configuration $Q\in\Qa$ is strongly connected iff its underlying heart configuration $H=H(Q)$ is. Note moreover that $H$ is distributed as in $\Ha$, by construction.

Recall the definition of s-cycle from the proof of Proposition~\ref{p:plain}, and note that if $H$ has no complex s-set of at most $m'/2$ arcs, then strong connectedness of $H$ is equivalent to $H$ having no s-cycles. Thus, in view of Lemma~\ref{l:heartnocomplex}, we only need to show that a heart configuration in $\Ha$ has no s-cycles with probability $1/9+o(1)$, and that when inserting $m-m'$ vertices in the arcs in order to generate a preheart configuration $Q\in\Qa$ we get a simple digraph a.a.s.

Since $\good$ holds, we have~\eqref{goodcond}. We first claim that this implies that a.a.s.\ the number $S$ of pairs of points that lie in the same in-bin is $O(r)$. Let $n_2:=a_{2,1}+a_{2,2}$ which must be $r-o(r)$. We have a distribution of $r+n_2$ points into $n_2$ in-bins chosen u.a.r.\ conditional upon  each bin receiving at least two points. If $r-n_2=o(\log r)$ say, immediately $S\le r + O(\log^2 r) = O(r)$. 
% query : the following part needs correlating with other similar arguments.
If on the other hand $r-n_2\to \infty$ (but recall it is $o(r)$), then this multinomial distribution can be approached by $n_2$ independent $2$-truncated Poissons conditional upon having sum $r+n_2$ (see Lemma~\ref{l:multinomial}). Combining Lemmas~\ref{l:variance} and~\ref{l:sum}, we deduce that $S=O(r)$ with probability $1-O((r-n_2)^{1/2}e^{-\log^3r})=1-o(1)$. 

The same holds for out-bins, so we may assume that the number of ways of choosing a set $\{a_1,a_2\}$ of out-points in the same bin and a set $\{b_1,b_2\}$ of in-points in the same bin is $O(r^2)$. The probability that $a_1,a_2$ are matched to $b_1,b_2$ thus creating a double arc is $O(1/r^2)$. The probability that a given double arc in $H$ gets no vertex inserted during the construction of $Q$ is $(m'-2)(m'-1)/(m-2)(m-1)=O(r^2/n^2)=o(1)$. Combining these conclusions, the expected number of double arcs in $H$ that get no vertex inserted during the construction of $Q$ is $o(1)$, and therefore $Q$ is simple a.a.s.

Let
\[
\mu_k = \sum_{j\ge1}^k \frac{2}{j} \left(\frac{2}{3}\right)^j
\quad\text{and}\quad
\mu = \lim_{k\to\infty}\mu_k = 2\log\frac{1}{1-2/3} = \log9.
\]
The number of s-cycles of order at most $k$ in $\Ha$ is asymptotically Poisson of mean $\mu_k$. This follows from estimating the factorial moments of this number of s-cycles in a similar way as in the proof of Proposition~\ref{p:plain}. The present case is simpler in two ways: firstly, there are no sets of vertices which are both a sink-cycle and a source-cycle, since this would imply having isolated cycles consisting of vertices of degree $(1,1)$. Secondly, the fact that the number of bins with degree exactly $(1,2)$ and the number of bins with degree exactly $(2,1)$ are each concentrated around $r$, and that the number of points in bins with higher degrees is negligible makes calculations much simpler than
those in the proof of Proposition~\ref{p:plain}. As before, the probability of having some s-cycles of order greater than $k$ can easily be bounded by some $f_k$ such that $\lim_{k\to\infty}f_k=0$. Therefore, the probability of having no s-cycles is $e^{-\mu}+o(1)$ as required.
\end{pf}
Finally, we proceed to prove Theorem~\ref{t:main} for the case $c\to1$. Denote by $K(n,m)$ the number of strongly connected digraphs with $n$ vertices and $m$ arcs. Given any degree sequence $\dv\in\hD$, there are exactly $m!(m'/m)$ preheart configurations in $\Qd$. Thus, in view of Lemma~\ref{l:heartuniform} and setting $A$ to be the event simple and strongly connected,  we can write
\begin{align}
K(n,m) &= \sum_{\dv\in\hD} \frac{m!(m'/m)\pr_\Qd(A)}{\prod_{i=1}^n d^+_i!d^-_i!}
\notag\\
&= (m-1)! \ex_{\hD}\left(m'\pr_\Qd(A)\right) \frac{(e^\lambda-1)^{2n}}{\lambda^{2m}} \pr_{\cD}\left(\Sigma\right)
\notag\\
& \sim \frac{(m-1)!}{2\pi (m-n)} \frac{(e^\lambda-1)^{2n}}{\lambda^{2m}} \ex_\hD\left(m'\pr_\Qd(A)\right),
\lab{eq:Knm}
\end{align}
since $\pr_\cD(\Sigma) \sim \frac{1}{2\pi nc(1+\lambda-c)} \sim \frac{1}{2\pi (m-n)}$ by Lemma~\ref{l:useful}.

To estimate $\ex_\hD\left(m'\pr_\Qd(A)\right)$ we will restrict ourselves to the event $\good$. If $\good$ holds, then~\eqref{goodcond} gives $m'\sim 3n'/2\sim 3(m-n)$. From Lemmata~\ref{l:heartnocomplex} and~\ref{l:heartsc}, for any any $\ba$ satisfying $\good$, we have 
\[
\pr_\Qa(A) \sim \frac19.
\]
Moreover, from Lemma~\ref{l:gooda}, we have that $\ex_\hD(m'(1-1_\good))=o(1)$ and in particular $\pr(\good)=1+o(1)$. Therefore,
\begin{align*}
\ex_\hD\left(m'\pr_\Qd(A)\right) &= \ex_\hD\left(m'1_\good\pr_\Qd(A)\right) + \ex_\hD\left(m'(1-1_\good)\pr_\Qd(A)\right)
\\
 &= (1+o(1))3(m-n) \pr_\Qnm(A \mid \good)
+ o(1)\\
&\sim (m-n)/3.
\end{align*}
Combining this with~\eqref{eq:Knm}, we obtain~\eqref{eq:main2} and thus complete the proof of the theorem.
%
%
%%%%%%%%%%%%%%%%%%%%%%%%%%%%%%%%%%%%%%%%%%%%%%%%%%
%
%
\section{Denser case: $c\to\infty$}\lab{s:dense}
In this section, we treat the case that $c\to\infty$ with $c=O(\log n)$. For such $c$, it follows easily from~\eqref{eq:ladef1} that
\bel{lambdacdense}
c=\lambda+o(1).
\ee
Our goal is to obtain the asymptotic number of strongly connected digraphs in this case, and therefore complete the proof of Theorem~\ref{t:main}. The main result in this section is the following.
\begin{prop}\lab{p:dense}
For $c:=m/n\to\infty$ with $c=O(\log n)$, a random digraph in $\dico$ is a.a.s.\ strongly connected.
\end{prop}
This result, combined with Theorem~\ref{t:dicoreenum}, gives the asymptotic number of strongly connected digraphs in the case that $c\to\infty$ with $c=O(\log n)$, which by~\eqref{lambdacdense} is asymptotic to~\eqref{eq:main3}, and thus the proof of Theorem~\ref{t:main} is complete.
%
%
% To show Prop number, change env to
% \begin{proof}[\bf Proof of Proposition~\ref{p:dense}]
%
\begin{pf}
As explained at the start of Section~\ref{s:sscd}, it suffices to show that a.a.s.\ there are no sink-sets. As before, we let $s$ be the cardinality of a hypothetical sink-set.
By duality it suffices to consider only $s\le n/2$.

Let $G\in  \dico$. We consider two cases. Let $K$ be fixed, and chosen sufficiently large as determined by the argument in Case~2 below.
\medskip

\noindent
{\em Case 1: $s\le c^K$}

Let $N_1$ be the number of vertices of outdegree $1$ in $G$. Define $f=\lambda /(e^\lambda-1)+n^{-1/3}$. The probability that one truncated Poisson r.v.\ equals 1 is $\lambda /(e^\lambda-1)$. Hence, by~\eqn{chern}, in the space  $\cD$ with $N_1$ interpreted in the natural way, we have
$$
\pr(N_1\ge f(c)n)=e^{-\Omega(\log^3 n)}.
$$
Then the same conclusion holds in $\hD$ by Lemma~\ref{l:useful}(a).
This also  transfers to the random graph $G\in  \dico$ by Lemma~\ref{l:useful}(b) and   Corollary~\ref{c:PtoG}, which show that probabilities multiply by  at most $e^{O(\log^2 n)}$.

We will condition on  $N_1=n_1$ where $n_1<f(c)n$, and consider the set $\Nscr=\{n_1:n_1<f(c)n,\; \pr(N_1=n_1)>1/n^2\}$. Then $\pr(N_1\notin \Nscr)<1/n=o(1)$. Let $H$ denote the event of not being strongly connected. We will show that
\[
\ \max_{n_1\in \Nscr} \pr(H\mid N_1=n_1) = o(1),
\]
which implies the result immediately. 
It helps to consider separately the event $J$ that the maximum in- or outdegree in the digraph is less than $\log^2 n$. By~\eqref{eq:tailbounds} and Lemma~\ref{l:sum}, $\pr(J)=1-o(1/n^2)$. 
So, letting $X_s$ denote the number of sink-sets of cardinality $s$, the above equation  follows if we show 
\bel{condprob2}
\sum_{1\le s\le n/2}\ex(X_s\wedge 1_{J}\mid N_1=n_1) = o(1)
\ee 
for all $n_1\in\Nscr$.

 By symmetry, we can assume the  $n_1$ vertices of outdegree 1 are specified in advance, so we may work in the restricted model, $\hatdico$, in which $V(G)$ is partitioned into two sets of vertices, $n_1$ in a set $A$ all of outdegree 1, and the rest in a set $B$  all of higher outdegrees. This is equivalent to $\dico$ conditioned on the set of vertices of outdegree 1 being precisely $A$. We use $\ex^*$ to denote expectation in this probability space.

 We will bound the probability $p(s,i,q)$  that $G\in D$ and some given set of vertices  $S$  is a sinkset of $G$ (also put $R=V\setminus S$), where $|S\cap A|=i$ and the set $Q$ of arcs with both ends in $S$ satisfies $|Q|=q$. The vertices in $B$ have outdegree at least 2, so $q \ge 2s-i$. By vertex symmetry,
\bel{estar}
\ex^* (X_s\wedge 1_J) =\sum_{i=0}^s{n_1\choose i}{n-n_1\choose s-i}\sum_{q\ge 2s-i}p(s,i,q)
\le n^s\sum_{i,q}(f(c))^ip(s,i,q).
\ee

We bound $p(s,i,q)$ using a switching technique. Take any digraph $G\in J$ with a sinkset $S$ as above, choose a set $Q'$ of arcs with both ends {\em not} in $S$ and with $|Q'|=q$, match up the arcs in $Q$ with those in $Q'$ in any manner, delete all arcs in $Q$ and $Q'$, and for each matched pair $uv \in Q$ and $u'v'$ in $Q'$, add the arcs $uv'$ and $u'v$. We call this operation a {\em switching}. The number of ways it can be performed on $G$ without creating any multiple arcs depends on the maximum degree (in- or out-) of $G$, which we denote by $\Delta$. For each arc $uv\in Q$, there are at most $\Delta$ arcs $wv$ and at most $\Delta^2$ arcs $wx$ excluded from choice as $u'v'$ due to causing double arcs. A similar number of exclusions of the form $wx$ come from arcs $ux$ also at least $m-2s\Delta$ arcs have both ends in $B$. Hence the number of valid switchings is at least $(m-O(\Delta^2+s\Delta))^q$. Performing such a switching produces some digraph $G'$ with the same (in,out)-degree sequence as $G$. How many such switchings can produce the same digraph $G'$? Assume $G$ has $r$ arcs directed from $R$ to $S$, so $G'$ has $r+q$ such arcs. Choose $q$ of these and pair them up with the  $q$ arcs of $G'$ from $S$ to $R$, to reverse the switching. This gives an upper bound (some reverse switchings  may be invalid) of say $(r+q)^q$ digraphs $G$ which can produce $G'$. 

At this point, we may deduce that the contribution to $p(s,i,q)$ from digraphs $G$ with the given values of $r$ and $\Delta$ is at most  
\bel{e:oneswitching}
\frac{(r+q)^q}{(m-O(\Delta^2+s\Delta))^q} = \left(\frac{O(1)(r+q)}{m}\right)^q,
\ee
since $G\in J$ and $\Delta<\log^2 n$. Note that the contribution to $p(s,i,q)$ from all $r$ such that $r\le c^3s^3$ or $r\le 18q$ is
\bel{e:oneswitching2}
\left(\frac{(c+q)^{O(1)}}{m}\right)^q.
\ee
To eliminate the influence of unusually large values of $r$ will require a more elaborate argument. If a vertex $v$ of $S$ in $G'$ is adjacent from $k$ vertices in $R$, where  $k>8c$, perform an additional switching to $G'$:   choose  $k- \lfloor 4c\rfloor$ arcs $u_1v,\ldots,u_{k-\lfloor 4c\rfloor}v$, $u_i\in R$, and replace them by arcs $u_iw_i$ with each $w_i\in R$ (without producing multiple arcs).  This produces a digraph $G_1$ having $A$ as its set of vertices of outdegree 1. The number of ways of performing this switching is at least (omitting floor functions for simplicity)
${k \choose k-4c }(n-s-\Delta)^{k-4c }$
since each vertex $u_i$ has outdegree at most $\Delta$. Each possible $G_1$ is produced in at most 
${m \choose k-4c }s$ ways, so the number of $G$ divided by the number of $G_1$ is at most 
$$\frac{{m \choose k-4c }  s }
{{k \choose k-4c}(n-s-\Delta)^{k-4c}}\le s(ec/k)^{k-4cs}(1+O(s/n+\Delta/n))^{ k-4cs}=O(s) (2ec/k)^{k/2},
$$
bounding the upper binomial above by $(em/(k-4c))^{k-4c}$ and the lower one below by $(k/(k-4c))^{k-4c}$, and using $8c<k\le \Delta<\log^2 n$.

If any vertex  of $G_1$ in $A$ is still adjacent from more than $8c$ vertices of $B$, we may repeat the previous step, to obtain $G_2$, and so on, up to some graph $G''$. If the switching is applied to    $j$ vertices of indegrees $k_1,\ldots, k_j$ with $\sum k_i = \tilde k$, the factors multiply to give
$$
s^j (ec)^{\tilde k / 2}\prod_i k_i^{-k_i/2}\le s^j (ec)^{\tilde k / 2}(j/\tilde k)^{-\tilde k/2}
$$
by   convexity of $x^x$. The worst case is $j=s$ which shows that  the number of possible $G'$ is   $s(O(1)c )/(\tilde k/s))^{\tilde k/2}$ times the number of $G''$, where $\tilde k\ge r-8cs$. Note that $G''\in \hatdico$. Thus, the contribution to $p(s,i,q)$ from $G$ with   given $r$   is at most
\bel{almost}
\frac{(O(1))^{q}(r+q)^q}{m^q}\left(\frac{O(1)cs}{r-8cs} \right)^{(r-8cs)/2}s^s,
\ee
where, by the conclusion using $G'$, the second factor can be taken to be $1$ for any particular $r$ (such as $r\le 8cs$). Note that $q\le s\Delta\le c^K\log^2 n$.  For $r$ larger than $c^3s^3$ and larger that $18q$, factor $(O(1)cs/(r-8cs))^{(r-8cs)/2}$ is at most $r^{-r/4}$ (for large enough $n$) and $s\le r^{r/18}$, so the product, summed over such $r$, is   $(o(1)/m)^q$.
In view of this and of~\eqref{e:oneswitching2},~\eqn{estar} gives
$$
\ex^* (X_s\wedge 1_J)  \le n^s\sum_{i,q}(f(c))^i((c+q)^{O(1)}/m)^q,
$$
where $i\le s$ and $q\ge 2s-i=s+s-i$. So this is at most
$$ 
(c^{O(1)}n/m)^s\sum_{i,t\ge 0}(f(c))^i(c^{O(1)}/m)^{s-i+t},
$$
which is at most $\big(c^{O(1)}\max\{f(c),1/m\}\big)^s$. Summing this over $s\le c^K$ gives $o(1)$ for the contribution to~\eqn{condprob2} from Case~1.

\medskip

%%%%%%%%%%%%%%%%%%%%%%%%%%%%%%%%%%%%%%%%%%
\noindent
{\em Case 2: $c^K< s\le n/2$}

Let $N^+_{\le 3}$ be the number of vertices of outdegree at most 3 in $G$, and let $h= c^3/(e^c-1)+c^3/n$. Then $\ex N^+_{\le 3}\le nh/2$, and by comparing with a binomial r.v.\ with expected value $nh/2$, and using~\eqn{chernsmallp}, we have
\[
\pr_{\cD}(N^+_{\le 3}\ge hn) = e^{-\Omega(hn)} = o(m^{-1}e^{-c^2}).
\]
So by Lemma~\ref{l:useful}(a,b), we deduce that $N^+_{\le 3}<hn$ a.a.s.\ in $\dico$, and it suffices to prove that,  conditional on this event, a.a.s.\ there are no sink-sets with cardinality $s$ in the range under consideration.

Henceforth, we consider $\Psnm$ (defined near the start of Section~\ref{s:sscd}) conditional upon a fixed outdegree sequence satisfying $N^+_{\le 3}=n^+_{\le 3}$ for some $n^+_{\le 3} <hn$. Again by Lemma~\ref{l:useful}(b), it is enough to show that if $R$ is the event that there exists a  sink-set  of size $s$ satisfying  $c^K\le s\le n/2$,
\bel{ProbR}
 \pr(R) =o(e^{-c^2}).
\ee
We note that our usual approach to proving properties of the degree sequence would be to work with independent truncated Poisson r.v.s for the degree sequences, prove what we want, and then condition on the sums. However, the last step increases probabilities of bad events in a manner unacceptable for the present argument. To avoid this, we define an auxiliary sequence $\hat d^-_1,\ldots,\hat d^-_n$ of independent copies of $\bindis(\hat n, \hat c/n)$, where $\hat n=(1+\delta)n$ and  $\hat c=(1+\delta)c$, and we set
\[
\delta=\epsilon/8,\ \epsilon=0.1.
\]
This sequence will be used to stochastically dominate some random variables defined on the indegree sequence of $\Psnm$. 

We next define some events that hold with high probability for the degree sequence.
 Let $\Delta = \lfloor5\log n/\log\log n+c^2\rfloor$ (so in particular $\Delta \le \log^3 n$, which suffices for most of our argument). This $\Delta$ turns out to be a typical bound on the maximum indegree. Let %
\begin{align*}
p_j&= \pr(\tpol=j)=\frac{\lambda^j}{(e^\lambda-1)j!},\\
\hat p_j&=\pr(\bindis(\hat n, \hat c/n)=j) = \binom{\hat n}{j}(\hat c/n)^j(1-\hat c/n)^{\hat n-j},
\end{align*}
and set
\bel{j0j3}
j_0=\min\{j\ge1\st n p_j\ge\log^{10}n\},\quad j_3=\max\{j \st n p_j\ge\log^{10}n\}.
\ee
Define the interval $I=\{j_0,\ldots,j_3\}$ and let $I'=\{1,\ldots,\Delta\}\setminus I$. Informally speaking, $I$ is the set of common indegrees and $I'$ the set of rare indegrees. Let $V$ and $V'$ be the set of vertices with indegrees in $I$ and $I'$, respectively. Define $H$ to be the event that the following holds: $d^-_{\max}\le\Delta$ (so $V$ and $V'$ partition the set of vertices); there exists a permutation $\sigma$ of $\{1,\ldots,n\}$ with the property that $d^-_i\le 1+\hat d^-_{\sigma(i)}$ for each $i\in V$; and moreover $|V'|=o(\log^{13}n)$. (The `$+1$' in the inequality $d^-_i\le 1+\hat d^-_{\sigma(i)}$ is to make it easier for our argument to cope with the fact that the Poisson variable is truncated at 1, whereas the binomial is not.)

We make several claims whose proofs are postponed. The first is the following.
\smallskip

\noindent{\bf Claim~1:} $\pr(H)=1-o(e^{-c^2})$.

The rest of the proof consists of showing that $\ex(X1_H)=O(0.93^s)$. This, together with Claim~1, gives~\eqn{ProbR} and we are done. 

%We now proceed to bound $\ex(X1_H)$.
 
Given a set $S$ of vertices ($|S|=s$), we may generate $\Psnm$ by specifying the random bijection $\phi$ last, which shows that 
\[
\pr(S\text{ is sinkset})\le \left(\frac{d^-(S)}{m}\right)^{d^+(S)}\le \left(\frac{d^-(S)}{m}\right)^{4s-3i}
\]
where $i$ is the number of vertices in $S$ of outdegree at most 3. Since the outdegree sequence was fixed and $n^+_{\le 3}<hn$, the number of sets $S$ of size $s$ with parameter $i$ is 
\[
\binom{n^+_{\le 3}}{i}\binom{n-n^+_{\le 3}}{s-i} \le \binom{s}{i}\frac{(hn)^in^{s-i}}{s!} \le 
\left(\frac{2en}{s}\right)^sh^i.
\]

Putting these together, we can bound the expected number $X$ of sink-sets of size $s$ restricted to the event $H$ by distinguishing cases according to the size of $d^-(S)$:
\bel{eq:key}
\ex(X1_H) \le \sum_{i=0}^s \left(\frac{2en}{s}\right)^sh^i \left[\left(\frac{(1+\epsilon)cs}{m}\right)^{4s-3i} + \sum_{t\ge(1+\epsilon)cs}^{\Delta s\wedge m} \left(\frac{t}{m}\right)^{4s-3i} \pr\big((d^-(S)=t)\wedge H\big)\right].
\ee

We need some care in treating the terms with $t\ge(1+\epsilon)cs$. Let $t_1=t(1+\epsilon/2) /(1+\epsilon)$ and $t_2=t\epsilon/\big(2(1+\epsilon)\big)$.  Note that $t_1+t_2=t$ with $t_1\ge(1+\epsilon/2)cs$ and $t_2\ge(\epsilon/2)cs$.

Let $S_I$ and $S_{I'}$ be the subsets in $S$ with degrees in $I$ and $I'$ respectively.
\smallskip

\noindent{\bf Claim~2:} $\pr\big((d^-(S_I)\ge t_1)\wedge H\big)\le e^{-\Omega(t)}$.

\smallskip

\noindent{\bf Claim~3:} $\pr\big((d^-(S_{I'})\ge t_2)\wedge H \big) = e^{-\Omega(t\log\log n)}$.%(s/(5t_2))^s$.
\smallskip
 
\noindent From these claims, it immediately follows that $\pr\big((d^-(S)=t)\wedge H\big)
= e^{-\Omega(t)}$,
and with~\eqref{eq:key} in mind we note that
%
%\[
%\ex(X1_H) \le \sum_{i=0}^s \left(\frac{2en}{s}\right)^sh^i \left[\left(\frac{(1+
% \epsilon)cs}{m}\right)^{4s-3i} + \sum_{t\ge(1+\epsilon)cs}^{\Delta s\wedge m} 
%\left(\frac{t}{m}\right)^{4s-3i} e^{-\Omega(t)}\right].
%\]
%
%Now, note that
\[
\sum_{t\ge(1+\epsilon)cs}^{\Delta s\wedge m} \left(\frac{t}{m}\right)^{4s-3i}e^{-\Omega(t)} = O\left(\left(\frac{(1+\epsilon)cs}{m}\right)^{4s-3i} e^{-\Omega(cs)}\right) = o\left(\left(\frac{(1+\epsilon)cs}{m}\right)^{4s-3i}\right),
\]
so then from~\eqn{eq:key} and using $m=cn$,
\bel{eq:key3}
\ex(X1_H) \le (1+o(1)) \sum_{i=0}^s \left(\frac{2en}{s}\right)^sh^i \left(\frac{(1+\epsilon) s}{n}\right)^{4s-3i}.
\ee
For $i\le s/100$ (recalling $\epsilon=0.1$ and $s\le n/2$),
\bean
\left(\frac{2en}{s}\right)^sh^i \left(\frac{(1+\epsilon)s }{n}\right)^{4s-3i}
&\le& \left( 2e (1+\eps) \right)^s  \left(\frac{(1+\epsilon)s }{n}\right)^{3s-3i}\\
&<& (2e(1.1))^s \left(\frac{1.1}{2}\right)^{2.97s} < 0.92^s,
\eean
whilst for $i\ge s/100$,
\[
\left(\frac{2en}{s}\right)^sh^i \left(\frac{(1+\epsilon) s}{n}\right)^{4s-3i}
\le \left(\frac{2en}{s}\right)^s h^{s/100} \left(\frac{1.1s}{n}\right)^s
\le \left(2.2e h^{1/100}\right)^s < 0.92^s.
\]
Thus,~\eqn{eq:key3} gives $\ex(X1_H)=O( s0.92^s)\le0.93^s$,
as desired. 
 It only remains to prove Claims~1--3.\smallskip
 
\noindent{\bf Proof of Claim~1.}
If $Y\distrib\tpol$, then
\[
\pr(Y\ge\Delta)=O(\pr(Y=\Delta)) = O((e\lambda/\Delta)^\Delta) = O(n^{-5}e^{-c^2}).
\]
Thus the statement holds for the first part of $H$ by taking union bound and using Lemma~\ref{l:useful}(a).

For the second part, consider a sequence $d^-_1,\ldots,d^-_n$ of independent truncated Poisson r.v.s, and a sequence $\hat d^-_1,\ldots,\hat d^-_n$ of independent Binomial r.v.s, as follows:
\[
d_i \distrib \tpol, \qquad \hat d_i \distrib \bindis(\hat n, \hat c/n),
\]
where $\hat n=(1+\delta)n$ and  $\hat c=(1+\delta)c$ (recalling $\delta=\epsilon/8$).

Recalling the definitions of $p_j$ and $\hat p_j$ above~\eqn{j0j3}, we have that
\[
\pr(d^-_i=j) = p_j \sim e^{-\lambda}\frac{\lambda^j}{j!}
\quad\text{and}\quad
\pr(\hat d^-_i=j) = \hat p_j \sim e^{-(1+\delta)^2c}\frac{((1+\delta)^2c)^j}{j!}.
\]
Let $Y_j$ and $\hat Y_j$ be the numbers of vertices of degree at least $j$ for each of the models, and similarly,   $Z_j$ and $\hat Z_j$  the numbers of vertices of degree at most $j$. We have
\[
\ex Y_j=n\pr(d^-_i\ge j), \quad
\ex \hat Y_j=n\pr(\hat d^-_i\ge j), \quad
\ex Z_j=n\pr(d^-_i\le j), \quad
\ex \hat Z_j=n\pr(\hat d^-_i\le j).
\]
Recall the definition of $j_0$ and $j_3$ in~\eqn{j0j3}, and let $j_1=c-\sqrt c/100$ and $j_2=(1+3\delta/2)c$. It is straightforward to check that $1\le j_0\le j_1\le j_2\le j_3\le \Delta$. If $j_2\le j\le j_3$, we easily verify that $\pr(d^-_i\ge j)=\Theta(p_j)$, and also that $p_j=o(\hat p_j)$ (by considering the ratios $p_{j+1}/p_j$ and $p_{j}/\hat p_j$). Hence, we have that $\pr(d^-_i\ge j)=o(\pr(\hat d^-_i\ge j))$.
If $j_1\le j\le j_2$, we have that $\pr(d^-_i\ge j)\le3/4$ (since $\mathrm{TPo}(\lambda)$ is asymptotically normal with mean $\lambda$ and variance $\lambda$, and truncation has a negligible effect on this), and clearly $\pr(\hat d^-_i\ge j)\sim1$ (for similar reasons or using second moment method). Therefore for 
$j_1\le j\le j_3$ we have that $\ex Y_j \le (4/5)\ex\hat Y_j$. Moreover, note that $Y_j$ and $\hat Y_j$ are binomially distributed, and for $j$ in this range we have that $\ex Y_j\ge np_j\ge\log^{10}n$. Hence, by~\eqref{chernsmallp} and taking a union bound, the probability that $|Y_j/\ex Y_j-1|>1/10$ or $|\hat Y_j/\ex\hat Y_j-1|>1/10$ for some $j$ in the range $j_1\le j\le j_3$ is $e^{-\Omega(\log^{10}n)}$. In particular, this implies that
$Y_j\le\hat Y_j$ for all $j\in [j_1, j_3]$ with probability $1-e^{-\Omega(\log^{10}n)}$.

On the other hand, if $j_0\le j\le j_1$, we easily verify that $\pr(\hat d^-_i\le j)=\Theta(\hat p_j)$, and also that $\hat p_j=o(p_j)$ (considering the ratios $\hat p_{j-1}/\hat p_j$ and $p_{j}/\hat p_j$). Therefore, $\ex\hat Z_j=o(\ex Z_j)$. Similarly as before, $Z_j$ and $\hat Z_j$ are binomially distributed and $\ex Z_j\ge np_j\ge\log^{10}n$. Using~\eqref{chernsmallp} again, we conclude that $\hat Z_j\le Z_j$ for all $j\in [j_0, j_1]$ with probability $1-e^{-\Omega(\log^{10}n)}$. Here we distinguish the two cases $\ex\hat Z_j\ge (\log^{10}n)/2$ and $\ex\hat Z_j<(\log^{10}n)/2$, and for the second case use stochastic domination of $\hat Z_j$ by a binomial r.v.\ of expectation $(\log^{10}n)/2$.

Summarising, we have that
\bel{almostthere}
Y_j\le\hat Y_j, \qquad\forall j\in [j_0+1, j_3]
\ee
with probability $1-e^{-\Omega(\log^{10}n)}$, where we used that $Z_j+Y_{j+1}=\hat Z_j+\hat Y_{j+1}=n$. However, what we really want is a suitable modification of~\eqref{almostthere} that holds for the range $j\in [j_0, j_3]$ and incorporates the `$+1$' shift in the definition of $H$. To do this, we distinguish two cases. If $j_0>1$, then it is straightforward to verify that $p_{j_0-1}n\ge\log^{8}n$, so the same argument as before but changing $j_0$ to $j_0-1$ shows that $Y_j\le\hat Y_{j}\le\hat Y_{j-1}$ for all $j\in [j_0, j_3]$ with probability $1-e^{-\Omega(\log^8n)}$. Otherwise if $j_0=1$, we trivially deduce from~\eqref{almostthere} that $Y_j\le\hat Y_{j-1}$ for all $j\in [j_0, j_3]$ with probability $1-e^{-\Omega(\log^{10}n)}$. Putting everything together, we conclude that
\[
Y_j\le\hat Y_{j-1}, \qquad\forall j\in [j_0, j_3]
\]
with probability $1-e^{-\Omega(\log^8n)}$. Thus, if this last inequality holds, then we can rearrange $\{1,\ldots,n\}$ by some permutation $\sigma$ in such a way that $d^-_i\le 1+\hat d^-_{\sigma(i)}$ for all $i\in I$.

Conditioning on the truncated Poisson r.v.s of the sequence $d^-_1,\ldots,d^-_n$ having fixed sum $m$ only multiplies the probability of failure by $O(\sqrt m)=o(n)$. The claim follows immediately.
\smallskip
 
\noindent{\bf Proof of Claim~2.}
Since we are restricting the probability space to the event $H$ and the choice of $S$ is uniformly randomised, we can bound the probability in question by replacing $\pr\big((d^-(S_I)\ge t_1)\wedge H\big)$ by $\pr(\hat d^-(S)+s\ge t_1)$, where $\hat d^-(S)=\sum_{i\in S}\hat d^-_i$ (informally speaking, $\hat d^-(S)+s$ is the total indegree of $S$ after having replaced the original indegree sequence by $\hat d^-_1+1,\ldots,\hat d^-_n+1$.)
In this model $\hat d^-(S)\distrib\bindis(s\hat n,\hat c/n)$, and it is immediate to verify (using standard deviation bounds on binomials; see also~\eqref{chernsmallp}) that $\pr(\hat d^-(S)\ge t_1-s)\le e^{-\Omega(t)}$, since $t_1/\ex \hat d^-(S)\ge(1+\epsilon/2)/(1+\delta)^2>1$ and $s=o(t_1)$.
\smallskip
 
\noindent{\bf Proof of Claim~3.}
Recall that $V'$ is the set of vertices with degrees in $I'$, and that $H$ implies that $|V'|\le\log^{13}n$. So the contribution of these to $d^-(S)$ is at most $\log^{16}n$. Thus, if $s>\log^{16}n$, we have $\pr(d^-(S_{I'})\ge t_2)=0$, since $t_2\ge\epsilon cs/2>\log^{16}n$.

So we may assume that $s\le \log^{16}n$. In this case $c\le\log^{16/K}n\le\log^{1/5}n$ (if say $K\ge 100$).  Thus, $c^2$ is negligible in the definition of $\Delta$ and we have $\Delta \sim 5\log n/\log\log n$. (Here we must be precise, since the trivial bound $\Delta\le\log^3n$ is not enough for this part of the argument.) Let us fix any $r\le\log^{13}n$ and restrict to the event that $|V'|=r$.  
Then we may use a model in which elements of $S$ are chosen independently, each with probability $s/n$, and condition  on  the size $s$ being achieved. Before conditioning, the number $Z=|V'\cap S|$ of red vertices in $S$ is $Z\distrib\bindis(|V'|,s/n)$, and conditioning on the size $s$ multiplies any probabilities by $O(s^{-1/2})$.  Note that $\ex Z=o(1)$ and thus, by elementary consideration of the binomial distribution, $\pr(Z\ge j) = O(\pr(Z=j))$.
Hence
\begin{multline*}
\pr(d^-(V'\cap S)\ge t_2\wedge H) \le \pr(Z\ge t_2/\Delta) = O(\pr(Z=t_2/\Delta)) =
\\
= O\left(\binom{|V'|}{t_2/\Delta}(s/n)^{t_2/\Delta}\right) =
O\left(\left(\frac{e|V'|s\Delta}{t_2n}\right)^{t_2/\Delta}\right) = O\left(e^{-\Omega(t\log\log n)}\right).
\qedhere
\end{multline*}
%\end{proof}
\end{pf}
%
%
%
%%%%%%%%%%%%%%%%%%%%%%%%%%%%%%%%%%%%%%%%%%%%%%%%%%%%%%
%
%
\section{Loop-free case}\lab{s:loopfree}
This section treats the case that digraphs are not permitted to have loops. We prove Theorems~\ref{t:mainloopfree} and~\ref{t:kdicoreenumloopfree}, which are analogues of Theorems~\ref{t:main} and~\ref{t:kdicoreenum}. To prove these theorems, we need the following result, which is similar to Lemmas~\ref{l:variance} and~\ref{l:concentration}.
\begin{lem}\lab{l:concentration2}
Let $Y^+_1,\ldots,Y^+_N,Y^-_1,\ldots,Y^-_N$ be independent r.v.s with $\tpo_k(\lambda)$ distribution, for fixed $k$ and for $0<\lambda\le\log N$. Let $c=\ex Y^+_1$. Then for any $t\ge\sqrt N\log^3N$ we have
\[
\pr\left(\Big|\sum_{i=1}^N Y^+_iY^-_i - c^2N\Big|>t\right) = O\left( e^{-(t^2/8N)^{1/5}} \right),
\]
asymptotically as $N\to\infty$.
\end{lem}
\begin{pf}
The argument is almost identical to that of Lemma~\ref{l:concentration}, so we just state the main differences. Here, we redefine $\Delta=(t^2/8N)^{1/5}$, $Y_{\max}=\max_{1\le i\le N} \{Y^+_i,Y^-_i\}$, $W_i = Y^+_iY^-_i-c^2$ and $W^*_i = W_i \,1_{E_i}$, where $E_i$ is the event that $Y^+_i\le\Delta$ and $Y^-_i\le\Delta$. Note that $\Delta=\Omega(\log^{6/5}N)$, and that $\lambda\le c\le(1+o(1))\log N$. It only remains to find appropriate bounds on $\pr(Y_{\max}>\Delta)$, $|\ex W^*_i|$ and $|W^*_i-\ex W^*_i|$, and then apply the same steps as in the proof of Lemma~\ref{l:concentration}.
The bound $\pr(Y_{\max}>\Delta) = O(e^{-\Delta})$ is obtained analogously. In view of~\eqref{eq:tailbounds} and the fact that $\ex((c-Y^-_i) \,1_{Y^-_i>\Delta})<0$, we easily deduce
\[
|\ex W^*_i| = \ex\left((c+Y^+_i) \,1_{Y^+_i\le\Delta}\right) \ex\left((Y^-_i-c) \,1_{Y^-_i>\Delta}\right)
\le 2c \ex\left(Y^-_i \,1_{Y^-_i>\Delta}\right) = O(e^{-\Delta}).
\]
Finally, we have $k^2-c^2 \le W^*_i \le \Delta^2-c^2$, and therefore $|W^*_i-\ex W^*_i|<\Delta^2$.
\end{pf}

\begin{proof}[\bf Proof of Theorem~\ref{t:kdicoreenumloopfree}]
 After extending Lemma~\ref{l:useful} to  the loop-free case, the proof is identical to that of Theorem~\ref{t:kdicoreenum}. So we just describe this extension of Lemma~\ref{l:useful}, which requires inserting an $e^{-c}$ factor in the asymptotic expressions in parts (b) and (c). The main adjustment in the proof is to redefine $\tilde F= \exp(-D_0-D^+D^-/2)$, where $D_0=\frac1m\sum_{i=1}^nd^+_id^-_i$. Instead of Lemma~\ref{l:PSimple},
we use a version which excludes loops. Again, we can use~\cite[Theorem~4.6]{M84} with digraphs  loop-free digraphs interpreted as bipartite graphs with a specific perfect matching being forbidden. Under the same conditions as Lemma~\ref{l:PSimple}, this implies that the probability that a random element of $\Pd$ has no loops and no multiple arcs is
\[
\exp\left(-\frac{1}{m}\sum_{i=1}^nd^+_id^-_i -\frac{1}{2m^2}\sum_{i,j=1}^nd^+_i(d^+_i-1)d^-_j(d^-_j-1) + O\left(\frac{\Delta^4}{m}\right)\right),
\]
uniformly for all $\dv$.
  $\Bad_2$ is the event that $|D_0-c|>t$ or $|D^+D^-/2-\lambda^+\lambda^-/2|>t$. We bound the probability that $|D_0-c|>t$ using Lemma~\ref{l:concentration2}.
\end{proof}
\begin{proof}[\bf Proof of Theorem~\ref{t:mainloopfree}.]
Again, we only need to point out how to change the proof of Theorem~\ref{t:main}.
For the case that $c$ is bounded and bounded away from $1$, we simply extend Propositions~\ref{p:complex} and~\ref{p:plain} in Section~\ref{s:sscd} to the loop-free case, and combine them with Theorem~\ref{t:kdicoreenumloopfree}. Proposition~\ref{p:complex} implies its own extension in this new setting, since the probability of an s-set when conditioning on loop-free digraphs can only increase by the inverse of the probability of having no loops, which is $\Theta(1)$ by comparing Theorems~\ref{t:kdicoreenumloopfree} and~\ref{t:kdicoreenum}.
Proposition~\ref{p:plain} is extended as follows.
\begin{prop}\lab{p:plainloopfree} Suppose that $c=m/n$ is bounded and bounded away from $1$. The probability
that a digraph in $\dico$
has no plain s-set is asymptotic to
\begin{equation}\lab{nononcomplex2}
e^{c(2/e^\lambda-1/e^{2\lambda})}
\frac{e^\lambda(e^\lambda-1-\lambda)^2}
{(e^{2\lambda}-e^\lambda-\lambda)(e^\lambda-1)},
\end{equation}
with $\lambda$  determined by the equation
$c=\lambda e^\lambda/(e^\lambda-1)$.
\end{prop}
\begin{pf}
The argument is almost identical to that of Proposition~\ref{p:plain}. We sketch the main differences. $C_k$ is again the number of s-cycles of order at most $k$ but we exclude s-cycles of order $1$ since they will be regarded as loops. $D$ is redefined to be the number of loops and double arcs. We have
\[
\mu_k = \sum_{j=2}^k \frac{2(c/e^\lambda)^j-(c/e^{2\lambda})^j}{j},
\]
and
\[
\mu=\lim_{k\to\infty}\mu_k = \log\left(\frac{(e^{2\lambda}-e^\lambda-\lambda)(e^\lambda-1)} {e^\lambda(e^\lambda-1-\lambda)^2}\right) - c(2/e^\lambda-1/e^{2\lambda}).
\]
The rest of the argument consists in bounding the probability of having s-cycles of order greater than $k$ for large $k$ and showing that $C_k$ and $D$ are asymptotically jointly independent Poisson with expectations $\ex_\Pnm C_k \sim \mu_k$ and $\ex_\Pnm D \sim c+\lambda^2/2$.
\end{pf}
The formula~\eqref{eq:main2}, for the very sparse case ($c\to1$) of Theorem~\ref{t:main}, remains unchanged: in the proof of Lemma~\ref{l:heartsc} one can easily see that the expected number of loops that get no vertices inserted while creating the preheart from the heart is $o(1)$ using an approach similar to the one for double arcs.

Finally, for the denser case ($c\to\infty$ with $c=O(\log n)$) it suffices to verify that Proposition~\ref{p:dense} is still valid if  loops are not permitted. Actually, the  argument in Section~\ref{s:dense} works for this setting with only the following trivial modifications. Note that for the first case in the proof ($s\le c^K$), the initial switchings do not create or destroy loops. The additional switchings can be performed in at least $\binom{k}{k-4c}(n-s-\Delta-1)^{k-4c}$ ways without creating loops (which only requires replacing $\Delta$ by $\Delta+1$) and the resulting bounds obtained in Section~\ref{s:dense} are unaffected.  The argument for the second case of the proof ($c^K<s\le n/2$) remains valid with the only difference that we have to additionally condition on having no loops. The extra effect of forbidding loops gives an additional   asymptotic $e^{-c}$ factor to the probability in Lemma~\ref{l:useful}(b) (see the proof of Theorem~\ref{t:kdicoreenumloopfree} for the extension of Lemma~\ref{l:useful} to loop-free digraphs).  Since $e^{-c^2}=o(e^{-c-\la^2/2})$, showing~\eqn{ProbR} still suffices in the loop-free context.
\end{proof}

\noindent
{\bf Acknowledgement\ } This paper relies heavily on results obtained in~\cite{PW1,PW2} by   Boris Pittel and the second author. We gratefully acknowledge discussions in 2001 with Boris on the subject of the present paper.

\end{document}